\documentclass[nohypdvips]{siamart251216}
\usepackage{amssymb,mathtools}
\usepackage[mathscr]{euscript}
\usepackage{tikz}
\definecolor{PaperInk}{HTML}{213747}
\definecolor{PaperAccent}{HTML}{2B6570}
\definecolor{PaperMuted}{HTML}{667882}
\definecolor{PaperRule}{HTML}{CBD7DC}
\newsiamremark{remark}{Remark}
\headers{Periodic gravity waves with a density jump}{Y. Fan, Y. Mao, and Q. Wang}
\makeatletter
\DeclareRobustCommand{\fint}{\mathpalette\paper@meanbar\relax\!\int}
\newcommand{\paper@meanbar}[2]{%
  \setbox\z@=\hbox{$\m@th#1\int$}%
  \vcenter{\hbox to\wd\z@{\hss$\m@th#1-$\hss}}%
  \kern-\wd\z@}
\makeatother
\NewDocumentCommand{\qty}{d() d[]}{%
\IfNoValueTF{#1}{%
\IfNoValueTF{#2}{\PackageError{qty}{Missing delimiter}{}}{\left[#2\right]}}%
{\left(#1\right)}}
\begin{document}
\title{Bifurcation and Nodal Structure of Periodic Gravity Waves with a Density Jump}
\author{Yanlong Fan\thanks{School of Mathematics, Southwest Jiaotong University, Chengdu, Sichuan, 611756, China.}
\and Yiqiu Mao\thanks{School of Mathematics and Information Science, Guangzhou University, Guangzhou, 510000, China.}
\and Quan Wang\thanks{Corresponding author. College of Mathematics, Sichuan University, Chengdu, Sichuan, 610065, China (\texttt{xihujunzi@scu.edu.cn}).}}
\maketitle

\begin{abstract}
We study the global bifurcation of periodic gravity waves in two stably stratified layers with fixed relative fluxes and a common Bernoulli offset. We first use the mean squared relative speed at the surface to remove the horizontal zero mode, including at critical points of the laminar Bernoulli function. Each transverse positive fundamental root then gives a local analytic branch. For nonpositive vorticities, a rational estimate gives uniqueness of the positive fundamental root when the upper-to-lower density ratio is at least \(3.05\times10^{-6}\). For each of two other parameter sets, interval estimates give at least three such roots. We next show that a laminar return on the bilateral component can occur at wave number \(k\ge1\), while a return on the crest-to-trough component must occur at a positive fundamental root. If this root is unique, the two orientations meet only at the base flow and together form the bilateral component. Finally, the depth estimate \(gM^2H_{\max}>1\), where \(M=\max h_p\), yields compactness on every bounded sublevel of \(M\) in the closure of the nonlaminar solutions. Thus, if the component through the bifurcation point contains no other admissible laminar flow, it contains a sequence along which \(\min(c-u)=1/M\to0\).
\end{abstract}
\begin{keywords}periodic water waves, density jump, transmission problem, global bifurcation, horizontal stagnation\end{keywords}
\begin{MSCcodes}35Q31, 35R35, 76B15\end{MSCcodes}

\section{Problem and principal conclusions}\label{sec_problem_results}
We consider two immiscible fluids of different densities bounded above by a free surface. A periodic gravity wave displaces both the surface and the interface. When the period and the relative flux in each layer are fixed, the mean depths may vary. The pressure is continuous across the interface, although the tangential velocity may have different traces there. Our aim is to determine how a connected family of periodic waves can continue beyond the small-amplitude regime. In particular, can the family lose compactness while the horizontal relative speed remains bounded away from zero? If the family returns to a laminar flow, which vertical mode is selected?

In the case of homogeneous rotational water waves, Constantin and Strauss \cite{ConstantinStrauss2004} obtained a global continuum by tracking monotonicity between the crest and the trough. Walsh \cite{Walsh2009} established global existence for continuously stratified periodic waves. For discontinuous density, Chen and Walsh \cite{ChenWalsh2016} obtained small periodic waves with piecewise smooth densities by perturbing the streamline density of a smooth wave. Walsh, B\"uhler and Shatah \cite{WalshBuhlerShatah2013} constructed small interfacial wind waves, including a finite atmosphere under a rigid lid; Matioc \cite{Matioc2012} treated internal waves with a critical layer. Wheeler \cite{Wheeler2019} obtained local periodic bifurcation with both a free surface and a density jump, allowing critical layers and imposing mean-depth and velocity constraints.

For large-amplitude waves, the boundary conditions and the choice of fixed flow parameters play an important role. Ambrose, Strauss and Wright \cite{AmbroseStraussWright2016} proved a continuation result for periodic vortex sheets with positive surface tension. Harmonic resonances in the capillary setting also produce Wilton ripples \cite{AkersAmbrose2025}. Doak, Sewell, Matthies and Wheeler \cite{Doak2026} obtained global periodic families with a vorticity interface between rigid walls. For solitary waves, Amick and Turner \cite{AmickTurner1986} treated the rigid-lid problem, and Wang \cite{Wang2017} constructed small solitary and generalized solitary waves with a free surface. Chen, Walsh and Wheeler \cite{ChenWalshWheeler2018} obtained near-stagnation along large-amplitude smoothly stratified solitary waves; Sinambela \cite{Sinambela2021} proved the corresponding conclusion with a density interface and prescribed upstream flow.

To our knowledge, the following question has remained open for periodic pure-gravity waves with a free surface and a stable density jump. Suppose that the fluxes are fixed, the vorticities are distinct constants in the two layers, and the layers have a common Bernoulli offset. Does a transverse positive fundamental root give a global crest-to-trough component on which every laminar return is fundamental and loss of compactness can occur only as the horizontal relative speed approaches zero? We give an affirmative answer under the condition \(c-u>0\). In fact, we obtain compactness under a uniform positive lower bound on \(c-u\) throughout the closure of all nonlaminar solutions. If the positive fundamental root is unique, we also show that the two orientations meet only at the base flow and together form the bilateral component.

We briefly describe the main difficulties and the arguments used to overcome them. First, at a turning point of the laminar Bernoulli value, the horizontal mean contributes an additional mode. We eliminate this mode by fixing the mean squared relative speed at the surface. The dispersion matrix for the surface and interface traces, together with its energy form, then determines the positive mode and its crossing. Second, where a symmetry axis meets the interface, a Hopf estimate in one layer alone cannot control both conormal traces. We construct a local barrier across the interface to obtain a weighted lower bound for \(-h_q\). Projecting the normalized limiting slope onto the first sine mode then shows that a laminar return must have the principal vertical profile on the oriented component. Finally, there is no prescribed far-field depth in the periodic problem. The local two-phase speed estimate degenerates as the boundaries approach one another. We first compare disk averages to separate the boundaries when \(M=\max h_p\) and \(Q\) are bounded. Testing the slope equation with the slope gives \(gM^2H_{\max}>1\). Since large \(Q\) at bounded \(M\) would force \(H_{\max}\to0\), a bound on \(M\) alone suffices for compactness.
\subsection{The height formulation}\label{subsec_height_formulation}

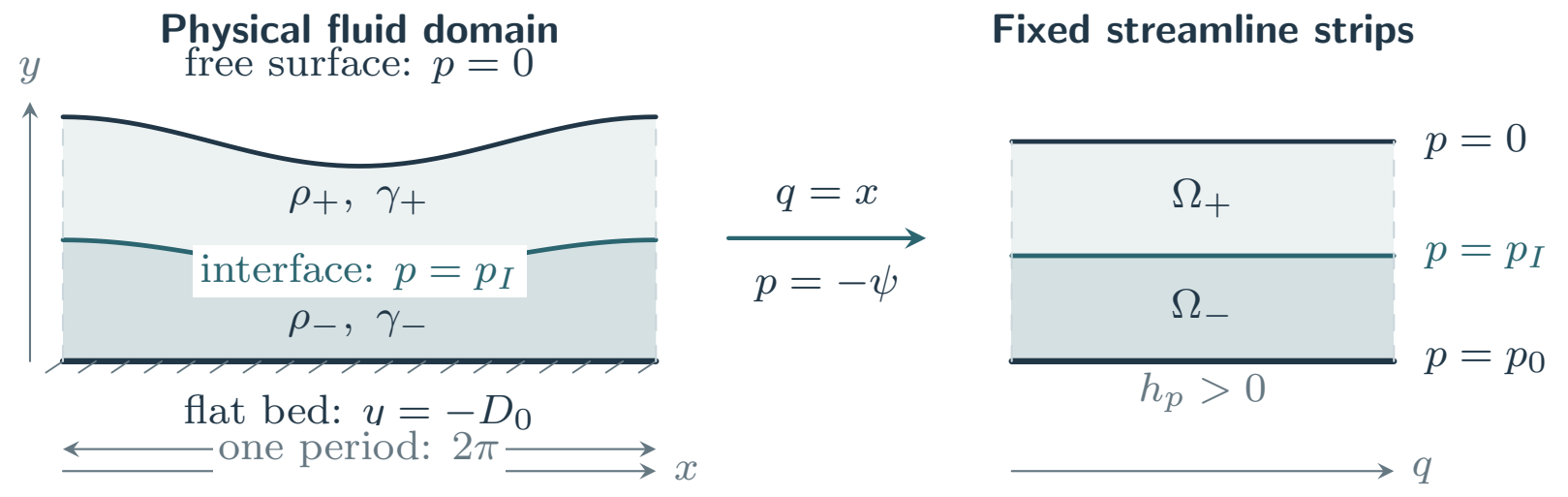
\begin{figure}[!ht]
\centering
\begin{tikzpicture}[x=.90cm,y=.60cm,>=stealth,font=\footnotesize,line cap=round,line join=round]
\node[font=\sffamily\footnotesize\bfseries,text=PaperInk] at (2.25,3.74) {Physical fluid domain};
\fill[PaperAccent!20] (0,0) -- (4.5,0) -- (4.5,1.38)
-- plot[domain=4.5:0,samples=55] (\x,{1.2+.18*cos(80*\x)}) -- cycle;
\fill[PaperAccent!9] plot[domain=0:4.5,samples=55] (\x,{1.2+.18*cos(80*\x)})
-- (4.5,2.78) -- plot[domain=4.5:0,samples=55] (\x,{2.5+.28*cos(80*\x)}) -- cycle;
\draw[PaperInk,line width=.9pt] plot[domain=0:4.5,samples=55] (\x,{2.5+.28*cos(80*\x)});
\draw[PaperAccent,line width=.9pt] plot[domain=0:4.5,samples=55] (\x,{1.2+.18*cos(80*\x)});
\draw[PaperInk,line width=1.2pt] (0,0)--(4.5,0);
\foreach \x in {0,.25,...,4.5} {\draw[PaperMuted,line width=.35pt] (\x,0)--(\x-.13,-.13);}
\draw[PaperRule,dashed] (0,0)--(0,2.78) (4.5,0)--(4.5,2.78);
\node[anchor=south,text=PaperInk] at (2.25,2.95) {free surface: \(p=0\)};
\node[text=PaperInk] at (2.25,1.83) {\({\rho}_{+},\ {\gamma}_{+}\)};
\node[text=PaperAccent,fill=white,inner sep=1.5pt] at (2.25,1.02) {interface: \(p={p_I}\)};
\node[text=PaperInk] at (2.25,.42) {\({\rho}_{-},\ {\gamma}_{-}\)};
\node[anchor=north,text=PaperInk] at (2.25,-.18) {flat bed: \(y=-{D_0}\)};
\draw[->,PaperMuted] (-.25,0)--(-.25,2.95) node[above] {\(y\)};
\draw[->,PaperMuted] (0,-1.25)--(4.5,-1.25) node[right] {\(x\)};
\draw[<->,PaperMuted] (0,-1.0)--(4.5,-1.0);
\node[fill=white,inner sep=1pt,text=PaperMuted] at (2.25,-1.0) {one period: \(2\pi\)};
\draw[->,PaperAccent,line width=.8pt] (5.05,1.4)--(6.55,1.4);
\node[above,text=PaperInk] at (5.8,1.5) {\(q=x\)};
\node[below,text=PaperInk] at (5.8,1.3) {\(p=-{\psi}\)};
\node[font=\sffamily\footnotesize\bfseries,text=PaperInk] at (8.65,3.74) {Fixed streamline strips};
\fill[PaperAccent!20] (7.2,0) rectangle (10.1,1.2);
\fill[PaperAccent!9] (7.2,1.2) rectangle (10.1,2.5);
\draw[PaperInk,line width=.9pt] (7.2,2.5)--(10.1,2.5);
\draw[PaperAccent,line width=.9pt] (7.2,1.2)--(10.1,1.2);
\draw[PaperInk,line width=1.2pt] (7.2,0)--(10.1,0);
\draw[PaperRule,dashed] (7.2,0)--(7.2,2.5) (10.1,0)--(10.1,2.5);
\node[text=PaperInk] at (8.65,1.85) {\({\Omega}_{+}\)};
\node[text=PaperInk] at (8.65,.6) {\({\Omega}_{-}\)};
\node[anchor=west,text=PaperInk] at (10.2,2.5) {\(p=0\)};
\node[anchor=west,text=PaperAccent] at (10.2,1.2) {\(p={p_I}\)};
\node[anchor=west,text=PaperInk] at (10.2,0) {\(p={p_0}\)};
\draw[->,PaperMuted] (7.2,-1.25)--(10.1,-1.25) node[right] {\(q\)};
\node[text=PaperMuted] at (8.65,-.35) {\({h}_p>0\)};
\end{tikzpicture}
\caption{One period in physical and streamline coordinates. The material interface separates the lighter upper fluid from the denser lower fluid, \(\rho_+<\rho_-\). Its pressure trace is continuous; its tangential velocity has two one-sided traces.}
\label{fig_two_layer_geometry}
\end{figure}

We take the negative streamfunction \(p=-\psi\) and the horizontal coordinate \(q=x\) as independent variables. Then \(h(q,p)\) measures height above the bed. The surface, interface and bed are the levels \(0,p_I,p_0\). The map in \cref{fig_two_layer_geometry} flattens the free surface and the interface simultaneously.
\begin{definition}\label{def_problem_data}
\leavevmode\par\smallskip\noindent
Fix \(g,D_0,c>0\), densities \(\rho_->\rho_+>0\), vorticities \(\gamma_-\ne\gamma_+\), and flux labels \(p_0<p_I<0\). Write \(I_-=(p_0,p_I)\), \(I_+=(p_I,0)\), \(\ell_i=|I_i|\), \(\Omega_i=\mathbb T\times I_i\), \(\Omega=\mathbb T\times(p_0,0)\), with \(\mathbb T=\mathbb R/(2\pi\mathbb Z)\). Set \(r=\rho_+/\rho_-\), \(\delta=1-r\) and \(\Delta\rho=\rho_+-\rho_-\).
Fix the H\"older exponent \({\alpha}\in\qty(0,1)\).
The subscript \(e\) denotes even periodic functions. All interface derivatives are one-sided traces, and \(\qty[f]_I=f\qty(\cdot,{p_I}+)-f\qty(\cdot,{p_I}-)\). Define
\begin{equation}\label{eq_height_spaces}
\begin{aligned}
{X}=\{u\in C^0\qty(\overline{{\Omega}}):\quad
&u_i\in C_e^{2,{\alpha}}\qty(\overline{{\Omega}_i}),\quad
u\qty(q,{p_0})=0\},\\
\|u\|_{{X}}&=\max_{i=\pm}
\|u_i\|_{C^{2,{\alpha}}\qty(\overline{{\Omega}_i})},\\
{Y}&=C_e^{0,{\alpha}}\qty(\overline{{\Omega}_{+}})
\times C_e^{0,{\alpha}}\qty(\overline{{\Omega}_{-}})
\times C_e^{1,{\alpha}}\qty({\mathbb T})^2,\\
{\mathcal U}&=\{u\in{X}:
\min_{i=\pm}\inf_{\overline{{\Omega}_i}}u_p>0\}.
\end{aligned}
\end{equation}
We equip the target space with the sum norm. For piecewise functions, \(\|u\|_{j,{\alpha};\mathrm{pw}}\) denotes the maximum of the two closed-layer \(C^{j,{\alpha}}\) norms.
\end{definition}

The unknowns are \(\qty({Q},{h})\in\mathbb R\times{\mathcal U}\). Write
\(b_+=\gamma_+p\), \(b_-=\gamma_-(p-p_I)+\gamma_+p_I\), \(B_i=Q/2-gD_0+b_i\), and \(\mathcal K(h)=(1+h_q^2)/(2h_p^2)\). The upper and lower interior residuals, surface residual and interface residual are
\begin{equation}\label{eq_model}
\mathcal F_i=(h_q/h_p)_q+(b_i-\mathcal K)_p,\quad
\mathcal F_S=(Q/2-gh-\mathcal K)_S,\quad
\mathcal F_I=[\rho(b-\mathcal K-gh+Q/2)]_I.
\end{equation}
Here \(\mathcal K=\mathcal K(h)\). The equation is \({\mathcal F}\qty({Q},{h})=0\). Its physical variables are
\begin{equation}\label{eq_physical_reconstruction}
\begin{aligned}
{\Phi}_{{h}}\qty(q,p)&=\qty(q,{h}\qty(q,p)-{D_0}),\qquad p=-{\psi},\\
{U}_i\circ{\Phi}_{{h}}
&=-{h}_p^{-1},\qquad
{V}_i\circ{\Phi}_{{h}}=-{h}_q/{h}_p,\\
{P}_i\circ{\Phi}_{{h}}
&={\rho}_i\qty({Q}/2+{b}_i\qty(p)
-{g}{h}-{\mathcal K}\qty({h})).
\end{aligned}
\end{equation}
The common Bernoulli offset varies with \({Q}\); the labels \({p_0},{p_I}\) remain fixed. Pressure continuity is imposed at the interface, while its tangential velocity traces may differ.

\subsection{Laminar parameters and the fundamental mode}\label{subsec_principal_laminar_data}

\begin{definition}\label{def_laminar_parameters}
Use the fluxes and density ratios of \cref{def_problem_data}.
For \({\lambda}>\max\qty(0,2{\gamma}_+{\ell}_+)\), define
\begin{equation}\label{eq_laminar_upper_depth_and_mu}
\begin{aligned}
{d}_+\qty({\lambda})
&=\frac{2{\ell}_+}
{\sqrt{{\lambda}}+\sqrt{{\lambda}-2{\gamma}_+{\ell}_+}},\\
{\mu}\qty({\lambda})
&={\lambda}-2{\gamma}_+{\ell}_+
+2{\delta}{g}{d}_+\qty({\lambda}).
\end{aligned}
\end{equation}
The laminar admissible set is
\begin{equation}\label{eq_laminar_admissible_set}
{\mathcal I}
=\left\{{\lambda}>\max\qty(0,2{\gamma}_+{\ell}_+):
{\mu}\qty({\lambda})>\max\qty(0,2{\gamma}_-{\ell}_-)\right\}.
\end{equation}
For \({\lambda}\in{\mathcal I}\), set
\begin{equation}\label{eq_laminar_lower_depth}
{d}_-\qty({\lambda})
=\frac{2{\ell}_-}
{\sqrt{{\mu}\qty({\lambda})}
+\sqrt{{\mu}\qty({\lambda})-2{\gamma}_-{\ell}_-}}.
\end{equation}
The associated laminar height and Bernoulli parameter are given by
\begin{equation}\label{eq_laminar_family}
\begin{aligned}
{H}_p\qty(p\mid{\lambda})
&=\begin{cases}
\qty({\mu}\qty({\lambda})
+2{\gamma}_-\qty(p-{p_I}))^{-1/2},
&p\in{I}_-,\\
\qty({\lambda}+2{\gamma}_+p)^{-1/2},
&p\in{I}_+,
\end{cases}\\
{H}\qty(p\mid{\lambda})
&=\int_{{p_0}}^p{H}_p\qty(t\mid{\lambda})\,dt,\\
{Q}\qty({\lambda})
&={\lambda}+2{g}
\qty({d}_-\qty({\lambda})+{d}_+\qty({\lambda})),\\
{H}_p\qty(0\mid{\lambda})^{-2}&={\lambda}.
\end{aligned}
\end{equation}
\end{definition}
\begin{proposition}\label{thm_laminar_family}
Equations \eqref{eq_laminar_upper_depth_and_mu}--\eqref{eq_laminar_family} give all admissible laminar solutions, once each. The parameter set is open and contains a right half-line. On each component the height is real analytic in \(\lambda\) with values in every closed-layer \(C^m\), and \(Q\) is analytic.
\end{proposition}
\begin{proof}
We set \(a=H_p^{-1}\). The interior equations give
\begin{equation}\label{eq_laminar_interior_constants}
 (a_i^2)'=2\gamma_i,\qquad a_S^2=\lambda,\qquad a_{I,-}^2=\mu .
\end{equation}
Writing \(U_S,U_I,V,V_B\) for the upper and lower endpoint speeds in descending order gives
\begin{equation}\label{eq_laminar_depth_integrals}
 d_i=\int_{I_i}a_i^{-1}\,dp=\frac{2\ell_i}{a_{i,\rm top}+a_{i,\rm bottom}} .
\end{equation}
The same formula holds when the vorticity vanishes. The free-surface condition gives \(Q=\lambda+2g(d_-+d_+)\). Pressure is hydrostatic in each layer, hence \(P_I=\rho_+gd_+\), and the lower Bernoulli law at the interface gives \(\mu=U_I^2+2\delta gd_+\). These identities give the claimed formulas. Conversely, direct substitution verifies the laminar equations. Both endpoint squared speeds must be positive in each layer, exactly the condition defining \(\mathcal I\). They tend to infinity with \(\lambda\). On compact admissible parameter neighborhoods the inverse square roots have uniformly convergent analytic expansions in every \(C^m\); integration preserves analyticity. Recovery of \(\lambda\) from the surface speed proves injectivity.
\end{proof}

\begin{definition}\label{def_modal_coefficients}
For \({\lambda}\in{\mathcal I}\), let \({a}_s\),
\({a}_{\mathrm I,+}\), and \({a}_{\mathrm I,-}\)
denote the positive relative speeds at the surface and on the two
sides of the interface.
Their values are
\begin{equation}\label{eq_modal_endpoint_speeds}
{a}_s=\sqrt{{\lambda}},\qquad
{a}_{\mathrm I,+}
=\sqrt{{\lambda}+2{\gamma}_+{p_I}},\qquad
{a}_{\mathrm I,-}=\sqrt{{\mu}}.
\end{equation}
For each integer \(k\ge1\), set
\begin{equation}\label{eq_modal_coefficients}
\begin{aligned}
{\mathsf A}_k
&=k{a}_s^2\coth\qty(k{d}_+)
-{\gamma}_+{a}_s-{g},\\
{\mathsf K}_k
&=k{a}_s{a}_{\mathrm I,+}
\operatorname{csch}\qty(k{d}_+),\\
{\mathsf D}_k
&=k{a}_{\mathrm I,+}^2\coth\qty(k{d}_+)
+{\gamma}_+{a}_{\mathrm I,+}\\
&\quad+\frac{{\rho}_-}{{\rho}_+}
\qty(k{a}_{\mathrm I,-}^2\coth\qty(k{d}_-)
-{\gamma}_-{a}_{\mathrm I,-})
+\frac{{\Delta\rho}}{{\rho}_+}{g},\\
{\mathsf M}_k
&=\begin{pmatrix}
{\mathsf A}_k&-{\mathsf K}_k\\
-{\mathsf K}_k&{\mathsf D}_k
\end{pmatrix},\qquad
{\mathscr D}_k=\det{\mathsf M}_k.
\end{aligned}
\end{equation}
\end{definition}

A parameter \(\lambda_*\in\mathcal I\) with
\(\mathscr D_1(\lambda_*)=0\) and \(\mathsf A_1(\lambda_*)>0\) is a
\emph{positive fundamental root}. The sign of \(\mathsf A_1\) selects a
profile that is positive above the bed. We use the following bifurcation
conditions
\begin{equation}\label{eq_bifurcation_conditions}
 \mathscr D_1(\lambda_*)=0,\qquad \mathsf A_1(\lambda_*)>0,\qquad
 \mathscr D_1'(\lambda_*)\ne0 .
\end{equation}
The derivative is taken along the laminar family. For a general height write
\begin{equation}\label{eq_gauge}
 \overline h_S=(2\pi)^{-1}\int_{\mathbb T}h(q,0)\,dq,\qquad
 \kappa=Q-2g\overline h_S .
\end{equation}
On a solution, \(\kappa\) is the surface mean of \((1+h_q^2)/h_p^2\); on a laminar solution it equals \(\lambda\).

\subsection{Principal results}\label{subsec_principal_results}

\begin{definition}\label{def_nodal_set}
In \(\mathscr O=\mathbb R\times\mathcal U\) let \(\mathscr Z=\mathcal F^{-1}(0)\), \(\Lambda\) its horizontal solutions, and \(\mathscr P\) its nonlaminar solutions with \(W=-h_q\ge0\) on \([0,\pi]\times[p_0,0]\). At the interface \(h_q\) is the derivative of the common trace.
For a base parameter \(\lambda_*\in\mathcal I\), write
\(z_*=(Q_*,H_*)=(Q(\lambda_*),H(\cdot\mid\lambda_*))\), and take the closures
relative to \(\mathscr O\) in the topology of \(\mathbb R\times X\):
\begin{equation}\label{eq_physical_solution_closures}
{\mathscr S}=\overline{{\mathscr Z}\setminus{\Lambda}}^{\,{\mathscr O}},
\qquad
{\mathscr S^+}=\overline{{\mathscr P}}^{\,{\mathscr O}}.
\end{equation}
When \(z_*\in{\mathscr S}\), let \({\mathcal C}\) be the
connected component of \({\mathscr S}\) through \(z_*\).
When \(z_*\in{\mathscr S^+}\), let \({\mathcal C^+}\) be the
connected component of \({\mathscr S^+}\) through \(z_*\).
In the latter case \({\mathcal C^+}\subseteq{\mathcal C}\).
\end{definition}

\begin{definition}\label{def_physical_controls}
For an actual solution \(\qty(Q,h)\in\mathscr Z\), use the physical
coordinate map \(\Phi_h\qty(q,p)=\qty(q,h\qty(q,p)-D_0)\) and streamfunction
\(\psi=-p\circ\Phi_h^{-1}\). Let \(D_\pm\) be the two physical layers
on the cylinder, bounded by the free surface \(S\), interface \(I\), and
bed \(B\). In each layer \(D_i\), let \(\psi_i=\psi|_{D_i}\) and write
\(u_i\) for the inertial-frame horizontal velocity. Define
\begin{equation}\label{eq_physical_control_quantities}
M=\max_{i=\pm}\|h_p\|_{L^\infty\qty(\overline\Omega_i)},\qquad
H_{\max}=\max_{q\in\mathbb T}h\qty(q,0),\qquad
V_I=\max_I|\nabla\psi_+|.
\end{equation}
The interface speed \(V_I\) uses the upper trace. By
\eqref{eq_physical_reconstruction}, \(M^{-1}\) is the minimum horizontal
relative speed on the two closed physical layers.
Define the vertical thicknesses \({d}_i\qty(q)\), the minimum boundary distance
\({d_*}\), and the interface slope bound
\({B_I}\) by
\begin{equation}\label{eq_layer_separation_definition}
\begin{aligned}
{d}_+\qty(q)&={h}\qty(q,0)-{h}\qty(q,{p_I}),
&{d}_-\qty(q)&={h}\qty(q,{p_I}),\\
{d_*}&=\min\{\operatorname{dist}\qty(I,S),\operatorname{dist}\qty(I,B)\},
&{B_I}&=\max_{q\in{\mathbb T}}|{h}_q\qty(q,{p_I})|.
\end{aligned}
\end{equation}
Distances are measured in the flat physical cylinder.
\end{definition}

\begin{theorem}\label{thm_local_bifurcation}
Suppose \eqref{eq_bifurcation_conditions} holds. There is a unique positive first-mode profile \(M_*\) with \(M_*(0)=1\). Set \(\varphi_*=M_*\cos q\) and \(\mathfrak b(w)=\pi^{-1}\int_{\mathbb T}w(q,0)\cos q\,dq\).
In a neighborhood of \(z_*=(Q(\lambda_*),H(\cdot\mid\lambda_*))\) the complete physical zero set consists of the laminar family and one analytic curve of even waves. Its nonlaminar points have unique parameter \(s=\mathfrak b(h(s))\), with \(0<|s|<\varepsilon\). Writing \(H_*=H(\cdot\mid\lambda_*)\), one has
 \(\kappa(s)=\lambda_*+O(s^2),\qquad h(s)=H_*+s\varphi_*+O_X(s^2),\qquad Q(s)=Q_*+O(s^2).\) 
The functions \(\kappa,Q\) are even in \(s\), and \(h(-s)(q,p)=h(s)(q+\pi,p)\).
\end{theorem}

\begin{theorem}\label{thm_bilateral_stagnation}
Suppose \(\lambda_*\in\mathcal I\) satisfies \eqref{eq_bifurcation_conditions}.
Then the base flow \(z_*\) belongs to \(\mathscr S\). On the bilateral
component \(\mathcal C\), each sublevel \(M\le K\), \(K>0\), is compact
in \(\mathbb R\times X\). At least one of the following holds.
\begin{enumerate}
\item For some \(\lambda_1\in\mathcal I\) with \(\lambda_1\ne\lambda_*\),
the flow \((Q(\lambda_1),H(\cdot\mid\lambda_1))\) belongs to \(\mathcal C\)
and \(\mathscr D_k(\lambda_1)=0\) for some \(k\ge1\).
\item There are \((Q_n,h_n)\in\mathcal C\) such that
\begin{equation}\label{eq_bilateral_horizontal_stagnation}
M(Q_n,h_n)\longrightarrow+\infty,\qquad
\min_{i=\pm}\min_{\overline D_{i,n}}(c-u_{i,n})
=\frac1{M(Q_n,h_n)}\longrightarrow0 .
\end{equation}
\end{enumerate}
\end{theorem}

\begin{theorem}\label{thm_oriented_stagnation}
Suppose \(\lambda_*\in\mathcal I\) satisfies \eqref{eq_bifurcation_conditions}.
The base flow \(z_*\) belongs to \(\mathscr S^+\).
Each nonlaminar member of \(\mathcal C^+\) satisfies
\begin{equation}\label{eq_oriented_strict_decrease}
h_q(q,p)<0\qquad(0<q<\pi,\ p_0<p\le0).
\end{equation}
This gives one crest and one trough per period on each streamline above
the bed, including the interface and free surface.
For each \(K>0\), the subset of \(\mathcal C^+\) with \(M\le K\) is
compact in \(\mathbb R\times X\). There are \(d_K>0\) and \(C_K<\infty\),
depending only on \(K,\alpha\) and the fixed data, such that on this subset
\begin{equation}\label{eq_oriented_compact_controls}
d_*\ge d_K,\qquad
Q+B_I+V_I+\max_{i=\pm}\|\nabla\psi_i\|_{L^\infty(\overline D_i)}
\le C_K .
\end{equation}
At least one of the following holds.
\begin{enumerate}
\item The component contains a different laminar flow at
\(\lambda_1\ne\lambda_*\) with \(\mathscr D_1(\lambda_1)=0\) and
\(\mathsf A_1(\lambda_1)>0\).
\item There is a sequence \((Q_n,h_n)\) on the component such that
\begin{equation}\label{eq_oriented_horizontal_stagnation}
M(Q_n,h_n)\longrightarrow+\infty,\qquad
\min_{i=\pm}\min_{\overline D_{i,n}}(c-u_{i,n})
=\frac1{M(Q_n,h_n)}\longrightarrow0 .
\end{equation}
\end{enumerate}
If the positive fundamental root is unique in the complete admissible
set \(\mathcal I\), then \eqref{eq_oriented_horizontal_stagnation} holds.
\end{theorem}

If the positive fundamental root is unique in \(\mathcal I\), then
 \(\mathcal C=\mathcal C^+\cup T_\pi\mathcal C^+,\qquad \mathcal C^+\cap T_\pi\mathcal C^+=\{z_*\}.\) 
Each orientation therefore contains a sequence approaching horizontal stagnation. For \(\gamma_\pm\le0\), the sufficient density ratio is \(r\ge3.05\times10^{-6}\) (\cref{prop_transversal_fundamental_regimes}). The two data sets in \cref{prop_fundamental_nonuniqueness} give three distinct local bifurcations in each case, one of them at a downward crossing.

\begin{figure}[htbp]
\centering\includegraphics[width=.90\linewidth]{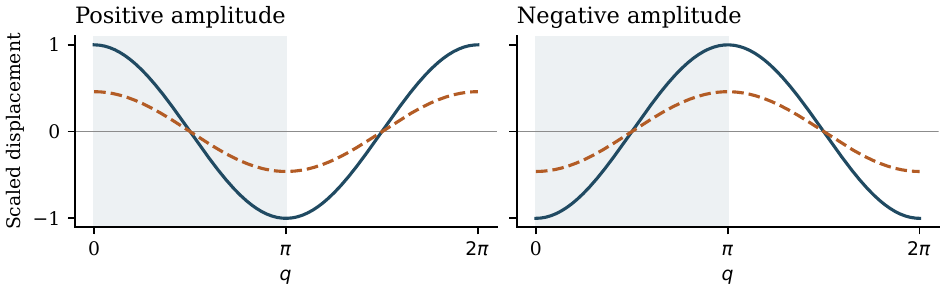}
\caption{First-order displacements from laminar height at the surface (solid) and interface (dashed), proportional to \(sM_*(p)\cos q\). The positive trace amplitudes are schematic and scaled for visibility. On the shaded half-period both traces decrease for positive amplitude and increase for negative amplitude. Translation by \(\pi\) exchanges the orientations.}
\label{fig_orientations}
\end{figure}

The leading term of the local branch explains the two orientations in \cref{fig_orientations}.

\section{Laminar flows and spectral structure}\label{sec_laminar_spectrum}
We first determine which laminar flows admit a nontrivial periodic perturbation. The horizontal mean requires a separate calculation since the Bernoulli parameter may have a critical point along the laminar family.

\subsection{The horizontal mean and its degeneracy}\label{subsec_zero_mode}

\begin{lemma}\label{lem_zero_mode}
For \({\lambda}\in{\mathcal I}\), define
\begin{equation}\label{eq_laminar_trace_integrals}
\begin{aligned}
{J}_+\qty({\lambda})
&=\int_{{p_I}}^0
\qty({\lambda}+2{\gamma}_+p)^{-3/2}\,dp,\\
{J}_-\qty({\lambda})
&=\int_{{p_0}}^{{p_I}}
\qty({\mu}\qty({\lambda})
+2{\gamma}_-\qty(p-{p_I}))^{-3/2}\,dp.
\end{aligned}
\end{equation}
Then
\begin{equation}\label{eq_zero_mode_parameter_derivatives}
\begin{aligned}
{\mu}'&=1-{\delta}{g}{J}_+,\\
{Q}'
&=1-{g}\qty({J}_++{J}_-)
+{\delta}{g}^2
{J}_+{J}_-.
\end{aligned}
\end{equation}
Linearize the height problem at \({H}\qty(\cdot\mid{\lambda})\), holding
\({Q}={Q}\qty({\lambda})\) fixed.
Its kernel among functions independent of \(q\) is nontrivial if and only if
\({Q}'\qty({\lambda})=0\). In that case this zero-mode kernel is
one-dimensional and is spanned by \(\partial_{{\lambda}}{H}\).
If \({Q}'\ne0\), the zero-mode kernel is \(\{0\}\).
\end{lemma}
\begin{proof}
Differentiation under the integral signs gives \(d_+'=-J_+/2,\ d_-'=-\mu'J_-/2\), hence the two derivative formulas. For a perturbation independent of the horizontal variable, let \(j_i=a_i^3m_p\), constant in each layer. If the Bernoulli value is allowed to vary by \(\dot Q\), the mean solution is uniquely determined by its surface squared-speed variation \(\dot\lambda=-2j_+\): integrating the upper equation fixes the depth variation, the pressure jump fixes the lower squared speed variation \(\dot\mu=(1-\delta gJ_+)\dot\lambda\), and integration from the bed fixes \(m=\partial_\lambda H\,\dot\lambda\). The surface equation gives \(\dot Q=Q'\dot\lambda\). Setting \(\dot Q=0\) proves the claim since the surface derivative of \(H_p^{-2}\) is one.
\end{proof}

Thus \(Q'(\lambda_*)=0\) detects degeneracy of the mean at fixed \(Q\). For bifurcation in the first harmonic we require \(\nu_1(\lambda_*)=-1\) and \(\nu_1'(\lambda_*)\ne0\). The surface-mean formulation also covers the case in which both types of degeneracy occur at the same parameter.

\subsection{Fourier modes and the dispersion matrix}\label{subsec_modal_matrix}

\begin{definition}\label{def_spectral_objects}
For \({\lambda}\in{\mathcal I}\), put
\({a}\qty(p)={H}_p\qty(p\mid{\lambda})^{-1}\).
On the full periodic spaces obtained by dropping the evenness restrictions
on \({X}\) and \({Y}\) in \cref{def_problem_data},
the laminar linearization is
\begin{equation}\label{eq_laminar_linearized_operator}
\begin{aligned}
\qty({L}_{{\lambda}}m)_i
&={a}\,m_{qq}+\qty({a}^3m_p)_p,\\
\qty({L}_{{\lambda}}m)_S
&={a}\qty(0)^3m_p\qty(q,0)-{g}m\qty(q,0),\\
\qty({L}_{{\lambda}}m)_I
&={\rho}_+{a}\qty({p_I}+)^3m_p^+
-{\rho}_-{a}\qty({p_I}-)^3m_p^-
-{\Delta\rho}{g}m\qty(q,{p_I}).
\end{aligned}
\end{equation}
All interface derivatives are one-sided traces. The coefficients are
independent of \(q\), so the expression preserves each parity class.
For local bifurcation we use its restriction
\({L}_{{\lambda}}:{X}\to{Y}\).
Let
\begin{equation}\label{eq_vertical_test_space}
{V_1}
=\{M\in H^1\qty({p_0},0):M\qty({p_0})=0\},
\end{equation}
and define, for \(M\in{V_1}\),
\begin{equation}\label{eq_modal_energy_definition}
\begin{aligned}
{E}_k\qty(M)
&=\sum_{i=\pm}{\rho}_i\int_{{I}_i}
\qty({a}^3M_p^2+k^2{a}M^2)\,dp\\
&\quad-{\rho}_+{g}M\qty(0)^2
+{\Delta\rho}{g}M\qty({p_I})^2.
\end{aligned}
\end{equation}
\end{definition}

\begin{proposition}\label{prop_mode_matrix}
For every \({\lambda}\in{\mathcal I}\) and integer \(k\ge1\),
the kernel of \({L}_{{\lambda}}\) in the even Fourier mode
\(M\qty(p)\cos\qty(kq)\) is nontrivial if and only if
\({\mathscr D}_k\qty({\lambda})=0\). Every such modal kernel is
one-dimensional. Its profile can be chosen strictly positive on
\(({p_0},0]\) if and only if \({\mathsf A}_k>0\).
If \({\mathsf A}_k<0\), every nonzero kernel profile has exactly one zero
in the upper layer and none in the interior of the lower layer.
\end{proposition}
\begin{proof}
For the mode \(M\cos kq\), \((a^3M_p)_p=k^2aM\). In physical height \(t=H(p)\), \(a_t=\gamma_i\). Thus \(z=aM\) satisfies \(z_{tt}=k^2z\) and \(a^3M_p=az_t-\gamma_i z\). On a layer of depth \(d\), measured from its bottom, the interpolant between values \(u,v\) is
\begin{equation}\label{eq_modal_interpolants}
 z(t)=\{u\sinh k(d-t)+v\sinh kt\}/\sinh kd.
\end{equation}
For \((\xi,\zeta)=(M(0),M(p_I))\) the two endpoint pairs from bottom to top are \((0,a_{\mathrm I,-}\zeta)\) and \((a_{\mathrm I,+}\zeta,a_s\xi)\). Substitution of the endpoint derivatives in the surface and interface residuals gives \(\mathsf A_k\xi-\mathsf K_k\zeta\) and \(\rho_+(\mathsf K_k\xi-\mathsf D_k\zeta)\). Since \(\mathsf K_k>0\), a singular matrix has rank one and \(\zeta/\xi=\mathsf A_k/\mathsf K_k\ne0\). If the traces have the same sign, normalization gives a strictly positive profile. Opposite signs give precisely one upper zero since \(z=C_1e^{kt}+C_2e^{-kt}\).
\end{proof}

\subsection{Energy, complete kernel, and transversality}\label{subsec_kernel_crossing}

\begin{lemma}\label{lem_mode_energy}
Fix \({\lambda}\in{\mathcal I}\), an integer \(k\ge1\),
and \(M\in{V_1}\). Put
\(\xi=M\qty(0)\), \(\zeta=M\qty({p_I})\), and, in the vertical
height coordinate \(t={H}\qty(p)\), let
\(z_i\qty(t)={a}\qty(p\qty(t))M\qty(p\qty(t))\), where \(p\qty(t)\) is the inverse
in the corresponding layer. Define \(z_{i,\mathrm h}\) by
\eqref{eq_modal_interpolants} with these endpoint values, and put
\(w_i=z_i-z_{i,\mathrm h}\). Then
\begin{equation}\label{eq_modal_energy_decomposition}
\begin{aligned}
{E}_k\qty(M)
&={\rho}_+
\begin{pmatrix}\xi&\zeta\end{pmatrix}
{\mathsf M}_k
\begin{pmatrix}\xi\\\zeta\end{pmatrix}\\
&\quad+{\rho}_-\int_0^{{d}_-}
\qty(w_{-,t}^2+k^2w_-^2)\,dt\\
&\quad+{\rho}_+
\int_{{d}_-}^{{d}_-+{d}_+}
\qty(w_{+,t}^2+k^2w_+^2)\,dt.
\end{aligned}
\end{equation}
\end{lemma}
\begin{proof}
Changing to \(t\) replaces the integrand by \(z_t^2+k^2z^2-\gamma_i(z^2/a)_t\). The remainder vanishes at the endpoints and is energy-orthogonal to the interpolant by \(z_{\mathrm h,tt}=k^2 z_{\mathrm h}\). The interpolant energy is \([z_{\mathrm h}z_{\mathrm h,t}]\); inserting its endpoint values, subtracting the vorticity boundary term and adding the gravitational traces gives exactly the matrix in \eqref{eq_modal_coefficients}. By continuity, the same calculation holds in \(H^1\) because the traces and the changes of coordinates in each layer are continuous.
\end{proof}

\begin{proposition}\label{prop_complete_kernel}
Suppose \({\lambda}_*\in{\mathcal I}\) satisfies
\begin{equation}\label{eq_positive_basic_root_conditions}
{\mathscr D}_1\qty({\lambda}_*)=0,
\qquad {\mathsf A}_1\qty({\lambda}_*)>0.
\end{equation}
Then \({L}_{{\lambda}_*}\) has no kernel in any integer
Fourier mode \(k\ge2\). If, in addition,
\({Q}'\qty({\lambda}_*)\ne0\), its complete kernel is
\begin{equation}\label{eq_complete_basic_kernel}
\begin{aligned}
\ker{L}_{{\lambda}_*}
&=\operatorname{span}\{{\varphi_*}\},\\
{\varphi_*}\qty(q,p)&={M_*}\qty(p)\cos q,\\
{M_*}\qty(0)&=1,\qquad
{M_*}\qty(p)>0\quad\qty({p_0}<p\le0).
\end{aligned}
\end{equation}
\end{proposition}
\begin{proof}
At the root the endpoint matrix is positive semidefinite of rank one, since \(\mathsf A_1>0\) and \(\mathsf D_1=\mathsf K_1^2/\mathsf A_1\). The interpolation energy therefore gives \(E_1(M)\ge0\). If \(k\ge2\) and \(M\ne0\),
\[
 E_k(M)=E_1(M)+(k^2-1)\sum_i\rho_i\int_{I_i}aM^2\,dp>0 .
\]
A kernel element has zero modal energy by integration by parts. At fixed \(Q\) the mean is excluded whenever \(Q'\ne0\) by \cref{lem_zero_mode}. The first-mode solution is positive, normalized and unique by \cref{prop_mode_matrix}; uniqueness of Fourier coefficients excludes every other kernel vector.
\end{proof}

\begin{proposition}\label{prop_transversality}
Under the hypotheses of \cref{prop_complete_kernel}, use its normalized
profile \({M_*}\) and vector \({\varphi_*}\).
The map \({\lambda}\mapsto{L}_{{\lambda}}\) is real analytic
with values in \(\mathcal L\qty({X},{Y})\), and
\begin{equation}\label{eq_modal_transversality_implication}
{\mathscr D}_1'\qty({\lambda}_*)\ne0
\quad\Longrightarrow\quad
{L}_{{\lambda}_*}'{\varphi_*}
\notin\operatorname{Range}{L}_{{\lambda}_*}.
\end{equation}
\end{proposition}
\begin{proof}
The laminar formulas are analytic in every layer coefficient norm on each admissible neighborhood. Fix \(\varphi_*=M_*\cos q\), \(v=(1,\mathsf A_1/\mathsf K_1)^T\) at the root, and define
\[
 \mathfrak L_* f=-\sum_i\rho_i\int_{\Omega_i}\varphi_*f_i
 +\rho_+\int_{\mathbb T}\varphi_{*,S} f_S-\int_{\mathbb T}\varphi_{*,I}f_I .
\]
Integration by parts gives the symmetric energy pairing, hence
\begin{equation}\label{eq_green_range_identity}
 \mathfrak L_*L_{\lambda_*}=0.
\end{equation}
When we hold the endpoint vector \(v\) fixed, the interpolant has energy \(\rho_+ v^T\mathsf M_1(\lambda)v\). At the root its profile variation pairs to zero with \(M_*\) by stationarity. Differentiating with respect to the coefficients therefore gives \(\mathfrak L_*L'_{\lambda_*}\varphi_*=\pi\rho_+v^T\mathsf M'_1v\), where \(\pi\) is the squared norm of the first harmonic. Since \(\operatorname{adj}\mathsf M_1=\mathsf D_1vv^T\) at the root, Jacobi's formula yields
\begin{equation}\label{eq_exact_transversality_pairing}
 \mathfrak L_*L'_{\lambda_*}\varphi_*=(\pi\rho_+/\mathsf D_1)\mathscr D'_1(\lambda_*).
\end{equation}
Here \(\mathsf D_1=\mathsf K_1^2/\mathsf A_1>0\), proving the range obstruction.
\end{proof}

\subsection{The principal vertical eigenvalue}\label{subsec_weighted_spectrum}

The weighted vertical eigenproblem associates the horizontal
wave numbers with signed spectral values. Its first eigenfunction
has a positive profile above the bed, so the first eigenvalue
at level \(-1\) identifies the positive fundamental roots.

\begin{definition}\label{def_full_spectrum}
For \({\lambda}\in{\mathcal I}\), use the laminar height
\({H}\qty(p\mid{\lambda})\) of
\cref{thm_laminar_family} and set
\begin{equation}\label{eq_full_spectrum_coefficients}
\begin{aligned}
{\mathfrak p}_i\qty(p,{\lambda})
&={\rho}_i{H}_p\qty(p\mid{\lambda})^{-3},&
{\mathfrak w}_i\qty(p,{\lambda})
&={\rho}_i{H}_p\qty(p\mid{\lambda})^{-1},\\
{\beta_S}&={\rho}_+{g}>0,&
{\beta_I}&={\Delta\rho}{g}<0.
\end{aligned}
\end{equation}
On the bed-zero space and the weighted Hilbert space
\begin{equation}\label{eq_full_spectrum_spaces}
\begin{aligned}
{V_1}
&=\left\{M\in H^1\qty({p_0},0):M\qty({p_0})=0\right\},\\
\mathcal H_{{\lambda}}
&=L^2\qty(\qty({p_0},0),{\mathfrak w}\,dp),\\
\qty(M,N)_{{\mathfrak w},{\lambda}}
&=\sum_{i=\pm}\int_{{I}_i}
{\mathfrak w}_i MN\,dp,
\end{aligned}
\end{equation}
define the symmetric form
\begin{equation}\label{eq_full_spectrum_form}
{\mathfrak a}_{{\lambda}}\qty(M,N)
=\sum_{i=\pm}\int_{{I}_i}
{\mathfrak p}_iM'N'\,dp
-{\beta_S}M\qty(0)N\qty(0)
+{\beta_I}M\qty({p_I})N\qty({p_I}).
\end{equation}
Its weak eigenvalue problem is
\begin{equation}\label{eq_full_spectrum_weak_problem}
M\in{V_1}\setminus\{0\},\qquad
{\mathfrak a}_{{\lambda}}\qty(M,N)
={\nu}\qty(M,N)_{{\mathfrak w},{\lambda}}
\quad\text{for every }N\in{V_1}.
\end{equation}
\end{definition}

\begin{proposition}\label{thm_full_spectrum}
The minimum of the Rayleigh quotient in \eqref{eq_full_spectrum_form} is a simple eigenvalue \(\nu_1\). Its eigenfunction can be chosen strictly positive above the bed. A nonnegative nonzero eigenfunction necessarily belongs to this eigenvalue. A positive fundamental root is equivalent to \(\nu_1=-1\).
\end{proposition}
\begin{proof}
On the bed-zero space the trace bound is \(|M(t)|^2\le\epsilon\|M'\|_2^2+C_\epsilon\|M\|_2^2\), by integrating from a nearby point and averaging. A minimizing sequence of weighted norm one is thus bounded in \(H^1\); compactness in \(C^0\) and lower semicontinuity give a minimizer. Its absolute value is also a minimizer. The Euler equation is \(-(\mathfrak p M')'=\nu\mathfrak w M\), with surface flux \(\beta_SM_S\) and flux jump \(\beta_I M_I\). At a zero of a nonnegative profile the flux vanishes: at the interface the two fluxes agree and the derivatives have opposite weak signs. At the surface, the same conclusion follows from the surface condition. Uniqueness for the ordinary differential equation then excludes every zero above the bed. The initial bed flux determines the profile up to scale. Orthogonality to the positive minimizer excludes other nonnegative eigenfunctions. The Fourier equation at \(k=1\) is this eigenproblem at \(-1\), and its trace ratio is \(\mathsf A_1/\mathsf K_1\).
\end{proof}
The principal eigenvalue varies continuously by the minimum principle and locally uniform coefficient comparison after a coercive shift. Since the principal eigenvalue is simple, it depends analytically on the parameter (analytic forms on a common domain, \cite{Kato1995}; the positive weight can equivalently be included on the right of the generalized eigenvalue equation). At a root normalize \(M_*(0)=1\). Differentiation of the weak eigenequation cancels the profile variations and gives
\begin{equation}\label{eq_fundamental_energy_derivative}
\begin{aligned}
\mathcal B_*={}&
\frac{\rho_+}{2}\int_{I_+}\qty(3a_+\qty(M_*')^2+a_+^{-1}M_*^2)\,dp\\
&+\frac{\rho_-\mu'\qty(\lambda_*)}{2}
\int_{I_-}\qty(3a_-\qty(M_*')^2+a_-^{-1}M_*^2)\,dp.
\end{aligned}
\end{equation}
The two derivatives are
\begin{equation}\label{eq_fundamental_root_derivatives}
\nu_1'\qty(\lambda_*)=
\frac{\mathcal B_*}{\displaystyle\sum_{i=\pm}\rho_i
\int_{I_i} a_iM_*^2\,dp},
\qquad
\mathscr D_1'\qty(\lambda_*)=\frac{\mathsf D_1\qty(\lambda_*)}{\rho_+}\mathcal B_*.
\end{equation}
To obtain the determinant identity, we differentiate the endpoint energy at fixed traces.
\begin{lemma}\label{thm_explicit_root_regime}
If \(\gamma_\pm\le0\), at least one positive fundamental root exists.
\end{lemma}
\begin{proof}
For \(\gamma_-\le0\) admissibility holds on the entire upper-layer half-line because \(\mu>0\). For two nonpositive vorticities this is \((0,\infty)\). At large \(\lambda\), all \(a_i^2\) lie between \(\lambda/2\) and \(2\lambda\). Since \(|M(t)|^2\le |p_0|\int |M'|^2\), both negative trace terms are absorbed by the derivative energy, and the Rayleigh quotient is bounded below by \(c\lambda\) for large \(\lambda\). Near zero take \(M=(1+p/\lambda)_+\). Its surface value is one, its derivative energy is \(O(\sqrt\lambda)\) because \(a_+^2\le(1+2|\gamma_+|)\lambda\) on its support, and its squared weighted norm is \(O(\lambda^{3/2})\). Thus \(\nu_1\to-\infty\) there. Continuity gives \(\nu_1=-1\).
\end{proof}
\subsection{A scalar residual and uniform uniqueness tests}\label{subsec_fundamental_root_tests}

\begin{proposition}\label{prop_fundamental_schur}
Suppose \(\gamma_-\le0\). Denote the endpoint speeds in
\eqref{eq_modal_endpoint_speeds} by \(U_S,U_I,V\), in order from
the surface to the lower trace at the interface, and set
\(V_B=\sqrt{\mu-2\gamma_-\ell_-}\). Write the upper endpoint
coefficients and the lower admittance as
\begin{equation}\label{eq_fundamental_admittances}
\begin{aligned}
N_{11}&=U_S^2\coth d_+-\gamma_+U_S,&
N_{12}&=-U_SU_I\operatorname{csch}d_+,\\
N_{22}&=U_I^2\coth d_++\gamma_+U_I,&
{n_-}&=\mu\coth d_--\gamma_-V.
\end{aligned}
\end{equation}
The entry \(\mathsf A_1=N_{11}-g\) is strictly increasing in
\(\lambda\). On \(\mathsf A_1>0\) define
\begin{equation}\label{eq_fundamental_schur_residual}
\begin{aligned}
{S_+}&=N_{22}-\frac{N_{12}^2}{\mathsf A_1},&
{F}&={n_-}+r{S_+}-\delta g,\\
{F}'&=r{S_+}'+\mu'{e_-},&
{e_-}&=\left.\partial_\mu{n_-}\right|_{\ell_-,\gamma_-}.
\end{aligned}
\end{equation}
Both \(S_+'\) and \(e_-\) are positive. On this set \(F\) has the
sign of \(\nu_1+1\), so its zeros are precisely the positive
fundamental roots. At each such root \(F'\), \(\mathscr D_1'\),
and \(\nu_1'\) have the same sign. All \(\lambda\) derivatives
hold the layer fluxes and the physical data fixed.
\end{proposition}
\begin{proof}
We eliminate the surface trace by minimizing the first-mode energy
at fixed interface trace. In the normalization of
\eqref{eq_modal_coefficients} the matrix is
\begin{equation*}
r\mathsf M_1=
\begin{pmatrix}
r\mathsf A_1&rN_{12}\\
rN_{12}&rN_{22}+n_--\delta g
\end{pmatrix},
\qquad
\mathscr D_1=\frac{\mathsf A_1F}{r}.
\end{equation*}
Completing the square in the surface trace leaves \(F\)
as the coefficient of the squared interface trace. The layer
interpolation in \cref{lem_mode_energy} splits off a nonnegative
energy with both traces zero, which vanishes only on the zero
function. Hence for \(F>0\) the first-mode energy is positive
definite, for \(F=0\) it is nonnegative with a nonzero null vector,
and for \(F<0\) it has a negative direction. The first Rayleigh
quotient is attained by \cref{thm_full_spectrum}. Its sign after
shifting the spectrum by \(1\) is therefore the sign of \(F\).

For the upper layer use \(s_+^2=\lambda+2\gamma_+p\) and the form
\begin{equation*}
\mathcal E_{+,\lambda}(m)=
\int_{I_+}(s_+^3m_p^2+s_+m^2)\,dp-gm(0)^2 .
\end{equation*}
The minimum at surface and interface traces \((1,0)\) is
\(\mathsf A_1\). For \(\mathsf A_1>0\), minimizing at interface
trace \(1\) with free surface trace gives \(S_+\). The analytic
minimizers are the layer interpolants at the minimizing traces.
In the derivative the variation of the profile vanishes by
stationarity, leaving in each case the positive integral
\begin{equation*}
\frac12\int_{I_+}(3s_+m_p^2+s_+^{-1}m^2)\,dp .
\end{equation*}
The lower admittance is the minimum of the lower integral energy
at traces \(0,1\) at the bed and interface. Write
\(s_-=H_p^{-1}|_{I_-}\) for the reciprocal laminar slope on this
layer. For the minimizer \(m_-\), differentiate at fixed
\(p,\ell_-,\gamma_-\) and then use
\(y=H(p)\) and \(dp=s_-\,dy\). Since \(\gamma_-\le0\) gives
\(s_-^2\ge\mu\), the lower identity and derivative satisfy
\begin{equation}\label{eq_lower_admittance_energy_bound}
\begin{aligned}
n_-&=\int_0^{d_-}s_-^2(m_{-,y}^2+m_-^2)\,dy,\\
0<e_-&=\frac12\int_0^{d_-}(3m_{-,y}^2+m_-^2)\,dy
\le\frac{3n_-}{2\mu}.
\end{aligned}
\end{equation}

At a positive fundamental root normalize its profile by \(m(0)=1\)
and write \(\zeta=m(p_I)=\mathsf A_1/\mathsf K_1>0\).
Minimization at this fixed interface trace gives the total
first-mode energy \(\rho_-\zeta^2F\). The energy
derivative \(\mathcal B_*\) in \eqref{eq_fundamental_energy_derivative}
and the differentiated determinant consequently satisfy
\begin{equation}\label{eq_fundamental_schur_derivative}
\begin{aligned}
\frac{\mathcal B_*}{\rho_-}
&=\zeta^2F'=rE_++\mu'\zeta^2e_-,
& E_+&=\frac12\int_{d_-}^{d_-+d_+}(3m_y^2+m^2)\,dy,\\
\mathscr D_1'&=\frac{\mathsf A_1}{r}F'.
\end{aligned}
\end{equation}
Equation \eqref{eq_fundamental_root_derivatives} now gives the
sign of \(\nu_1'\).
\end{proof}

\begin{proposition}\label{prop_transversal_fundamental_regimes}
Assume \(\gamma_-\le0\).
\begin{enumerate}
\item If
\begin{equation}\label{eq_nonpositive_vorticity_density_regime}
\gamma_+\le0,\qquad r\ge3.05\times10^{-6},
\end{equation}
then there is a unique positive fundamental root \(\lambda_*\) in
\(\mathcal I\), with \(\mathscr D_1'(\lambda_*)>0\) and \(\nu_1'(\lambda_*)>0\).
It satisfies \eqref{eq_bifurcation_conditions}.

\item For general \(\gamma_+\), put \(\gamma_+^+=\max\{\gamma_+,0\}\)
and \(\widehat d=\min\{5/4,\sqrt{2\ell_+/|\gamma_+|}\}\), with
\(\widehat d=5/4\) when \(\gamma_+=0\). Let \(\lambda_t\) be the unique
parameter with \(\delta gJ_+(\lambda_t)=1\). Each of the following
conditions gives at most one positive fundamental root, with
\(\mathscr D_1'>0\) and \(\nu_1'>0\) at any such root:
\begin{equation}\label{eq_general_upper_root_tests}
\begin{aligned}
rg&\ge(\gamma_+^+\sqrt{\widehat d}+\sqrt{\delta g})^2
            (\widehat d\coth\widehat d-1),\\
\mathsf A_1(\lambda_t)&\le0,\\
d_+(\lambda_t)&\ge5/4 .
\end{aligned}
\end{equation}
The second condition also gives existence. A sufficient condition for
the first is \(rg\ge 0.4736(1.1181\gamma_+^++\sqrt{\delta g})^2\).
\end{enumerate}
\end{proposition}
\begin{proof}
We bound the potentially negative lower-layer term in
\eqref{eq_fundamental_schur_derivative} by the upper-layer energy
and use the resulting strictly upward crossings of \(\nu_1+1\)
to count roots.

\smallskip\noindent\textit{Step 1. Bounds when the lower term is negative.}\enspace
In the notation of \eqref{eq_fundamental_schur_derivative}, we show
that at a positive fundamental root with \(\mu'<0\) and \(d_+\ge5/4\),
 \(rE_+>|\mu'|\zeta^2e_-.\) 
We shall also use these trace estimates when one of the layers is thin.

At a root with \(\mu'<0\) write \(d=d_+\), \(x=U_SU_I\) and
\(\varphi(d)=d\coth d-1\). Integration in \eqref{eq_laminar_trace_integrals}
gives \(J_+=d/x\). The inequalities \(\mu'<0\) and \(\mathsf A_1>0\)
and the identity \(U_I=U_S-\gamma_+d\) give
\begin{equation}\label{eq_negative_mu_geometry}
0<x<\delta gd,\qquad
d\mathsf A_1=U_S^2\varphi(d)+x-gd>0,\qquad
U_S^2\varphi(d)>rgd .
\end{equation}
The surface normalization and the interface trace \(\zeta\) give
\begin{equation}\label{eq_upper_fundamental_energy_lower_bounds}
E_+\ge\max\left\{\frac12\tanh d,\frac{3(1-\zeta)^2}{2d}\right\}.
\end{equation}
For the first bound use the minimizer \(\cosh y/\cosh d\) of
\(\int_0^d(m_y^2+m^2)\,dy\) at \(m(d)=1\). The second follows
from the Cauchy--Schwarz inequality applied to the difference of traces.

The interface row of the matrix gives
\begin{equation*}
\begin{aligned}
\frac{|\mu'|}{\mu}&<\frac1{2x},&
n_-&=\delta g-rN_{22}+r\mathsf K_1/\zeta,\\
\zeta^2n_-&\le\delta g\zeta^2+r\mathsf K_1\zeta,&
\mathsf K_1&=x\operatorname{csch}d,\qquad
\mathsf K_1\zeta=\mathsf A_1 .
\end{aligned}
\end{equation*}
Here \(\mu=U_I^2+2\delta gd\), and \(N_{22}>0\) is the value of the
upper integral energy with zero surface trace and unit interface
trace. For \(d>1/2\), put \(\Delta=U_I^2+\delta g(2d-1)>0\).
Since \(n_-\ge\mu\), the interface row gives
\(\zeta<r\mathsf K_1/\Delta\). Substitution in
\eqref{eq_lower_admittance_energy_bound} yields
\begin{equation}\label{eq_fundamental_thick_negative}
\begin{aligned}
|\mu'|\zeta^2e_-
&<\frac{3}{4x}(\delta g\zeta^2+r\mathsf K_1\zeta)\\
&<\frac34r^2\operatorname{csch}^2d
              \left(\frac{\delta gx}{\Delta^2}+\frac{x}{\Delta}\right)\\
&\le\frac34r^2\operatorname{csch}^2d\,m_d\left(1+\frac1{2d-1}\right).
\end{aligned}
\end{equation}
For general \(\gamma_+\) take \(m_d=d/(2d-1)\), since
\(x<\delta gd\). When \(\gamma_+\le0\) we also have \(x\le U_I^2\),
so take \(m_d=\min\{1,d/(2d-1)\}\). For \(d\ge5/4\) and any
\(\gamma_+\), the bound in \eqref{eq_fundamental_thick_negative}
is smaller than \(r\tanh d/2\) by
 \(3rd^2\operatorname{csch}^2d<(2d-1)^2\tanh d .\) 
At \(d=5/4,r=1\) this is \eqref{eq_certificate_cutoff_values}.
Afterwards \(d/\sinh d\) strictly decreases while the right side grows.

\smallskip\noindent\textit{Step 2. Nonpositive upper vorticity.}
For \(\mu'\ge0\), \eqref{eq_fundamental_schur_derivative} already gives \(F'>0\). Suppose \(\gamma_+\le0\) and \(\mu'<0\). Put
 \(u=gd/x-1,\quad \phi=d\coth d-1,\quad b=\sinh(d)/d,\quad a=b\phi,\qquad m=2(1+\phi)/(3+2\phi).\) 
Now \(U_S^2\le x\le U_I^2\), \(x<\delta gd\), and \(0<d\mathsf A_1\le x(\phi-u)\). Thus \(0<u<\phi\), \(\zeta\le b(\phi-u)\), \(\mathsf A_1<\delta g(\phi-u)\), and \(|\mu'|=\delta u-r\le u\).
Since \(\mu\ge x+2\delta gd\) and \(x=gd/(1+u)\), we have \(2\delta gd/\mu\le m\). Using \(\zeta^2n_-\le\delta g\zeta^2+r\mathsf A_1\) gives
\[
 \frac{|\mu'|\zeta^2 e_-}{r}
 \le\frac{3m}{4d}\left(\frac{b^2}{r}\,u(\phi-u)^2+u(\phi-u)\right)
 \le\frac{3m}{4d}\left(\frac{4b^2\phi^3}{27r}+\frac{\phi^2}{4}\right).
\]
Indeed \(4\phi^3-27u(\phi-u)^2=(4\phi-3u)(\phi-3u)^2\) and \(\phi^2-4u(\phi-u)=(\phi-2u)^2\). If \(a<1\), \eqref{eq_upper_fundamental_energy_lower_bounds} thus proves upward crossing whenever
\begin{equation}\label{eq_fundamental_thin_crossing_tests}
 T(d,r):=\frac{m}{2(1-a)^2}\left(\frac{4b^2\phi^3}{27r}+\frac{\phi^2}{4}\right)<1 .
\end{equation}
For the complementary comparison put \(L(d)=2d-1+\tanh d\). Since \(U_I^2\ge x>gd/(1+\phi)=g\tanh d\), the denominator \(\Delta=U_I^2+\delta g(2d-1)\) is greater than \(\delta g L(d)\). If \(L>0\), the same interface-row estimate \(\zeta<r\mathsf K_1/\Delta\) used in Step 1 applies even below \(d=1/2\). The first two lines of \eqref{eq_fundamental_thick_negative} and \(x<\delta gd\) yield
\[
 \frac{|\mu'|\zeta^2 e_-}{r}< B(d,r):=
 \tfrac34 r d\,\operatorname{csch}^2(d)(L^{-2}+L^{-1}).
\]
On a band \(l\le d\le h\), use \(E_+\ge\tanh(l)/2\), and also \(E_+\ge3(1-a(h))^2/(2h)\) if \(h<\infty\) and \(a(h)<1\). Write \(E(l,h)\) for the largest applicable bound. The functions \(\phi,b,a,m\) increase: the numerator of \(\phi'\) is \(\sinh d\cosh d-d>0\), and \(d\cosh d-\sinh d>0\) gives \(b'>0\). Also \(d\,\operatorname{csch}^2 d\) decreases because \(2d\coth d>1\), and \(L\) increases. Hence \(T\) is least favorable at the greatest depth and least density ratio, whereas \(B/E(l,h)\) is least favorable at the least depth and greatest ratio. The five finite bands in \cref{sec_fundamental_root_certificates} together with Step~1 for \(d\ge5/4\) prove strict crossing for all \(r\ge3.05\times10^{-6}\).
Existence follows from \cref{thm_explicit_root_regime}. Every zero of the continuous function \(\nu_1+1\) on the admissible half-line crosses upwards. Two such zeros would leave a positive point followed by a negative one; the first intervening zero could not cross upwards. This proves uniqueness.

\smallskip\noindent\textit{Step 3. General upper vorticity in part (2).}\enspace
We use each test in \eqref{eq_general_upper_root_tests} to exclude
thin-layer roots with a negative lower contribution. The second test
also gives existence.

For a root with \(\mu'<0\) and \(d<5/4\), solve
\(U_S^2-\gamma_+dU_S=x<\delta gd\) for \(U_S\). Then
\eqref{eq_negative_mu_geometry} gives
\begin{equation*}
rg<\frac{U_S^2\varphi(d)}d
 <(\gamma_+^+\sqrt d+\sqrt{\delta g})^2\varphi(d).
\end{equation*}
The last quantity is increasing in \(d\), and
\((U_S+U_I)^2\ge |U_S^2-U_I^2|=2|\gamma_+|\ell_+\) gives
\(d\le\widehat d\) in this case. Thus the first condition of
\eqref{eq_general_upper_root_tests} excludes these roots.
The simpler sufficient condition follows from
\eqref{eq_certificate_cutoff_values}.

The derivative formula gives \(\mu'=1-\delta gJ_+\). The integral defining \(J_+\) decreases strictly from infinity to zero on the upper admissible half-line, so \(\delta gJ_+(\lambda_t)=1\) has a unique solution. Hence \(\mu'<0\) exactly for
\(\lambda<\lambda_t\). If \(\mathsf A_1(\lambda_t)\le0\), the strict
increase of \(\mathsf A_1\) from \cref{prop_fundamental_schur}
excludes positive fundamental roots at \(\lambda\le\lambda_t\).
To see existence, at \(\lambda_t\)
take the trace vector \((1,0)\) if \(\mathsf A_1<0\), or
\((1,\zeta)\) for sufficiently small \(\zeta>0\) if
\(\mathsf A_1=0\). In the latter case the energy at the
interpolant is negative because \(\mathsf K_1>0\), hence
\(-2\mathsf K_1\zeta+\mathsf D_1\zeta^2<0\). In both cases
\(\nu_1(\lambda_t)<-1\), whereas \(\nu_1\to+\infty\) by the
large-parameter estimate in the proof of
\cref{thm_explicit_root_regime}.
For \(d_+(\lambda_t)\ge5/4\), the strict decrease of \(d_+\)
puts every negative-\(\mu'\) root in the thick-layer case.
The assertions now follow by counting the crossings as above.
\end{proof}

\begin{proposition}\label{prop_fundamental_nonuniqueness}
Fix \(g=1\). List the data as \((r,\ell_+,\ell_-,\gamma_+,\gamma_-)\).
For the data
 \((3/10,\,1/2,\,1/2,\,10,\,0)\) 
the following three disjoint intervals each contain exactly one
positive fundamental root:
\begin{equation}\label{eq_positive_shear_root_brackets}
[10.0000048,10.0000050],\quad
[10.00268,10.00270],\quad
[10.00565,10.00568].
\end{equation}
For the data
 \((1.5\times10^{-6},\,0.1913,\,3.0967,\,0,\,-10^{-8})\) 
the corresponding three intervals are
\begin{equation}\label{eq_nonpositive_shear_root_brackets}
[0.3270,0.3275],\qquad[0.329,0.3295],\qquad[0.331,0.3318].
\end{equation}
In each example the crossings are upward, downward, and upward,
in that order, and the admissible set is a half-line. Each data set therefore has at least three distinct bifurcating laminar flows; each has the analytic branch and global alternatives of \cref{thm_local_bifurcation}--\cref{thm_oriented_stagnation}.
\end{proposition}
\begin{proof}
The enclosing calculations in
\cref{sec_fundamental_root_certificates} give
\(\mathsf A_1>0\) on the six intervals, opposite signs for \(F\) at the
two endpoints of each, and the indicated sign of \(F'\) throughout each bracket. The intermediate value theorem and \cref{prop_fundamental_schur} give the assertions. Both sets of data satisfy
\(\gamma_-\le0\), so the laminar speed formulas give
the admissible half-lines.
\end{proof}

\begin{figure}[htbp]
\centering\includegraphics[width=.95\linewidth]{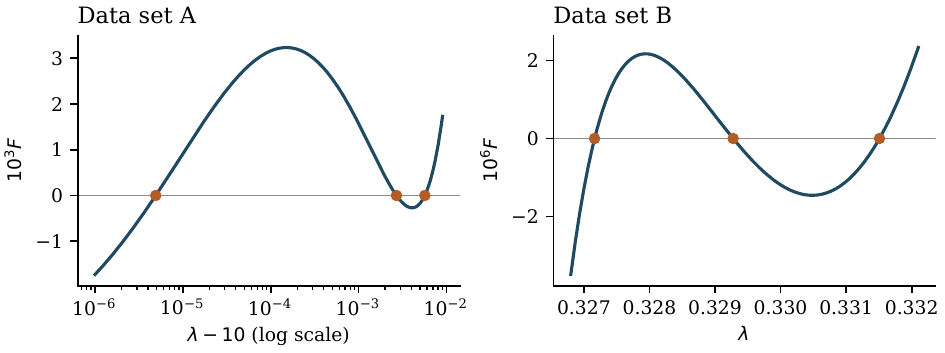}
\caption{Scalar residual \(F\) for the first (A) and second (B) data sets of \cref{prop_fundamental_nonuniqueness}. Dots mark the numerical roots in the three certified brackets for each example. In increasing order of \(\lambda\), each curve crosses zero upwards, downwards and upwards. The logarithmic horizontal scale in A resolves the first root near \(\lambda=10\).}
\label{fig_roots}
\end{figure}

Figure~\ref{fig_roots} illustrates the crossings for both sets of data in \cref{prop_fundamental_nonuniqueness}.

\subsection{Rational depth bounds and root enclosures}\label{sec_fundamental_root_certificates}
\paragraph{The uniqueness constant.}
We regard all terminating decimals in this calculation as exact rational numbers. We set \(P_j(x)=\sum_{k=0}^4 x^{2k}/(2k+j)!\) and \(U_j(x)=P_j(x)+2x^{10}/(10+j)!\), \(j=0,1\). For \(0<x\le5/4\), the ratio of consecutive omitted terms is at most \((5/4)^2/(11\cdot12)<1/2\). Thus \(P_0\le\cosh x\le U_0\) and \(P_1\le b(x)\le U_1\). The table gives rational bounds \(p\ge U_0/P_1-1\ge\phi(x)\), \(v\ge U_1\ge b(x)\), and \(t\le xP_1/U_0\le\tanh x\). Each row covers \((l,x]\), with \(l\) the preceding endpoint, initially zero.
\[
\begin{array}{c|ccc|c|cc}
x&p&v&t&r_0&T_{\rm up}&(B/E)_{\rm up}\\\hline
.3382&.03784&1.0192&.3258&.00000305&.999&--\\
.35&.04051&1.0206&.3363&.000008&.472&.934\\
.50&.08198&1.0422&.4621&.0001&.364&.068\\
.80&.20476&1.1102&.6640&.01&.100&.092\\
1.25&.47357&1.2816&.8482&.1&.761&.401
\end{array}
\]
We give the details for the last two columns. For \(r\ge r_0\), substitute \((p,v,r_0)\) for \((\phi,b,r)\) in \eqref{eq_fundamental_thin_crossing_tests}; \(vp<1\) in each row. In the first row, where the margin is smallest, the prefactor is less than \(.3651\), and the parenthesis less than \(2.7342\), whose product is below \(.999\). For \(r\le r_0\) in rows 2--5, use the preceding row's \(t_l\). We set \(D=2l-1+t_l>0\) and \(e=\max\{t_l/2,3(1-vp)^2/(2x)\}\). Since \(\sinh l\ge l\), the rational expression \(3r_0(D^{-2}+D^{-1})/(4le)\) bounds \(B(l,r_0)/E(l,x)\) from above. For example, in row 2, \(D=.0022\) and \(e>3.938\), giving a ratio below \(.934\). These substitutions require only multiplication and division of positive rationals. Finally \(P_1(5/4)>1.28\), and the last row gives, at \(d=5/4\),
\begin{equation}\label{eq_certificate_cutoff_values}
 3/b^2<3/1.28^2<(3/2)^2(.8482)\le(2d-1)^2\tanh d,
 \quad \phi<.4736,\quad \sqrt d<1.1181 .
\end{equation}
Step~1 covers all greater depths. The comparisons therefore establish the sufficient bound \(r\ge3.05\times10^{-6}\) without numerical root finding or interval software.

\paragraph{Root intervals.}
These separate multiple-root examples use outward-rounded interval arithmetic (\texttt{verify.py}, \texttt{mpmath} 1.3.0, 40 decimal digits). With \(L,R\) the endpoints of a bracket in
\eqref{eq_positive_shear_root_brackets} or
\eqref{eq_nonpositive_shear_root_brackets}, evaluate
\(\mathsf A_1(L),F(L),F(R)\) and enclose \(F'\) on the whole bracket.
Write \(\mathsf K_1=U_SU_I\operatorname{csch}d_+\).
Differentiating at fixed fluxes gives
\begin{equation}\label{eq_certificate_closed_form_derivatives}
\begin{aligned}
d_+'&=-\frac{d_+}{2U_SU_I},\\
\mathsf A_1'&=\coth d_+-U_S^2\operatorname{csch}^2d_+\,d_+'-\frac{\gamma_+}{2U_S},\\
\mathsf K_1'&=\mathsf K_1\left(\frac{1}{2U_S^2}+\frac{1}{2U_I^2}-d_+'\coth d_+\right),\\
N_{22}'&=\coth d_+-U_I^2\operatorname{csch}^2d_+\,d_+'+\frac{\gamma_+}{2U_I},\\
e_-&=\coth d_-+\frac{\mu d_-\operatorname{csch}^2d_-}{2VV_B}-\frac{\gamma_-}{2V},\\
F'&=\mu'e_-+r\left(N_{22}'-\frac{2\mathsf K_1\mathsf K_1'}{\mathsf A_1}
                        +\frac{\mathsf K_1^2\mathsf A_1'}{\mathsf A_1^2}\right).
\end{aligned}
\end{equation}
Each interval lies inside its admissible half-line, and
\(\mathsf A_1(L)>0.0136\) in all six intervals. Its strict increase
then gives \(\mathsf A_1>0\) throughout. The certificate gives
signs \(F(L),F(R),F'\) equal to \((- ,+,+)\), \((+,-,-)\),
and \((- ,+,+)\) in that order for each example. The endpoint
magnitudes exceed \(4.4\cdot10^{-7}\), and the enclosures for
\(F'\) on all subintervals give \(|F'|>4.3\cdot10^{-5}\).
For each rational bracket, the verification script continues bisection until the signed derivative has the claimed lower bound on every piece. Forty further sign bisections enclose the root for evaluation of \(Q'=1-gd_+/(U_SU_I)-g\mu'd_-/(VV_B)\).

At the three roots in each example the laminar derivative is enclosed, respectively, in
\[
\begin{array}{c|ccc}
 &1&2&3\\ \hline
 A&(7.84,7.85)&(-0.332,-0.331)&(-0.457,-0.455)\\
 B&(0.0470,0.0473)&(0.0262,0.0265)&(0.0047,0.0048).
\end{array}
\]

In case A the opposite signs of \(Q'\) at the first two roots give an intervening critical point. All six root enclosures exclude \(Q'=0\). Case B has \(r=1.5\times10^{-6}\); the optimal uniform density-ratio threshold is undetermined.

\section{The analytic transmission framework}\label{sec_analytic_framework}
In this section we establish the estimates for the linear transmission problem and then apply them to the nonlinear height operator. The estimates will also be used to pass to limits along the bifurcating component.

\subsection{Linear transmission estimates, solvability and Fredholm index}\label{subsec_linear_estimates}

\begin{definition}\label{def_transmission_system}
Use the intervals and strips of \cref{def_problem_data}. For this linear
system the two densities need only satisfy \({\rho}_\pm>0\), and
\({\sigma_S},{\sigma_I}\in\mathbb R\) are arbitrary fixed
Robin coefficients. Set
\begin{equation}\label{eq_transmission_spaces}
\begin{aligned}
{\widehat X}
&=\{u\in C^0\qty(\overline{{\Omega}}):
u_i\in C^{2,{\alpha}}\qty(\overline{{\Omega}_i}),\
u\qty(q,{p_0})=0\},\\
\|u\|_{j,{\alpha};\mathrm{pw}}
&=\max_{i=\pm}\|u_i\|_{C^{j,{\alpha}}\qty(\overline{{\Omega}_i})},\\
{\widehat Y}
&=C^{0,{\alpha}}\qty(\overline{{\Omega}_+})
\times C^{0,{\alpha}}\qty(\overline{{\Omega}_-})
\times C^{1,{\alpha}}\qty({\mathbb T})^2,\\
{V_2}
&=\{u\in H^1\qty({\Omega}):u\qty(q,{p_0})=0\}.
\end{aligned}
\end{equation}
The norm on \({\widehat X}\) is the displayed piecewise norm, and the norm
on \({\widehat Y}\) is the sum of its four component norms.
All functions are periodic in \(q\).

In the coordinate order \(\qty(q,p)\), let symmetric matrices
\({\mathbf C}_i\in
C^{1,{\alpha}}\qty(\overline{{\Omega}_i};\mathbb R^{2\times2})\)
satisfy
\begin{equation}\label{eq_transmission_coefficients}
\begin{aligned}
{\nu_0}|\xi|^2
&\le \xi^T{\mathbf C}_i\qty(q,p)\xi
\le {K}|\xi|^2,\\
\max_{i=\pm}\|{\mathbf C}_i\|_{C^{1,{\alpha}}}
&\le {K},
\qquad {\nu_0},{K}>0.
\end{aligned}
\end{equation}
For \({\tau}\ge0\), define
\begin{equation}\label{eq_transmission_operator}
\begin{aligned}
\qty({\mathscr T}_{{\mathbf C},{\tau}}u)_i
&=\operatorname{div}\qty({\mathbf C}_i\nabla u)
-{\tau}{\rho}_i u,\\
\qty({\mathscr T}_{{\mathbf C},{\tau}}u)_S
&=\qty({\mathbf C}_+\nabla u)_p\qty(q,0)
-{\sigma_S}u\qty(q,0),\\
\qty({\mathscr T}_{{\mathbf C},{\tau}}u)_I
&=\qty[\qty({\mathbf C}\nabla u)_p]_{{p_I}}
-{\sigma_I}u\qty(q,{p_I}).
\end{aligned}
\end{equation}
The jump is the upper trace minus the lower trace; all interface
derivatives are one-sided.
When needed, impose the reflection covariance
\begin{equation}\label{eq_transmission_reflection_covariance}
J_q=\operatorname{diag}\qty(-1,1),
\qquad
{\mathbf C}_i\qty(-q,p)=J_q{\mathbf C}_i\qty(q,p)J_q.
\end{equation}
Even and odd periodic subspaces of \({\widehat X}\) and \({\widehat Y}\)
are understood componentwise.

\end{definition}

\begin{lemma}\label{lem_transmission_gradient}
Use the flat strips and Robin coefficients of \cref{def_transmission_system}.
Suppose the symmetric matrices \({\mathbf C}_i\) are uniformly
elliptic with lower bound \({\nu_0}>0\) and have piecewise
\(C^{0,{\alpha}}\) norm bounded by \({K}\).
Let \(u\in{V_2}\) be piecewise \(C^{1,{\alpha}}\) and
satisfy the following equations in the weak sense.
\begin{equation}\label{eq_transmission_first_order_data}
\begin{aligned}
\operatorname{div}\qty({\mathbf C}_i\nabla u)
&=\operatorname{div}G_i+f_i,\\
\qty({\mathbf C}_+\nabla u-G_+)_p|_{p=0}
&={\sigma_S}u_S+g_S,\\
\qty[\qty({\mathbf C}\nabla u-G)_p]_{{p_I}}
&={\sigma_I}u_I+g_I,
\qquad u|_{p={p_0}}=0,
\end{aligned}
\end{equation}
where \(G_i\in C^{0,{\alpha}}\), \(f_i\in L^\infty\) and
\(g_S,g_I\in C^{0,{\alpha}}\qty({\mathbb T})\).
There is a constant \(C\) depending only on the strip geometry, the
Hölder exponent, the matrix bounds and the two Robin coefficients such that
\begin{equation}\label{eq_transmission_gradient_estimate}
\begin{aligned}
\|u\|_{1,{\alpha};\mathrm{pw}}
\le C\bigl(&\|u\|_{L^2}
+\|G\|_{0,{\alpha};\mathrm{pw}}+\|f\|_{L^\infty}\\
&+\|g_S\|_{C^{0,{\alpha}}}
+\|g_I\|_{C^{0,{\alpha}}}\bigr).
\end{aligned}
\end{equation}
\end{lemma}
\begin{proof}
We choose upper-layer cutoffs \(\chi_S,\chi_I\) equal to one at the indicated boundary and zero near the other. Extend them by zero below. With boundary loads \(b_j=\sigma_j u_j+g_j\), set
\begin{equation}\label{eq_transmission_flux_lifting}
 F_{\rm lift}=e_p(\chi_S b_S+\chi_I b_I),\quad J=G+F_{\rm lift},\quad f_0=f-\partial_p(F_{\rm lift})_p
\end{equation}
layerwise. The effective flux \(C\nabla u-J\) has zero jump and zero surface normal component, so its distributional divergence is \(f_0\). We use even reflection at the surface and odd reflection at the bed. Writing \(S=\operatorname{diag}(1,-1)\) and \(\epsilon=+1,-1\) respectively, the reflected fields are \((\epsilon u,SCS,\epsilon SJ,\epsilon f_0)\). Since the normal fluxes agree and the bed trace vanishes, both extensions satisfy the reflected equation in the weak sense. The jumps in the normal variable do not affect these tangential H\"older estimates.

We cover the strips by finitely many balls of fixed radius. On this cover, Dong's laminate estimate \cite[Corollary~2.7]{Dong2012} bounds \(\|\nabla u\|_\infty\) and the \(C^\alpha\) seminorms of \(u_q\) and \(U=(C\nabla u-J)_p\) by \(C(\|u\|_2+\|J\|_\infty+[J]_{q;\alpha}+\|f_0\|_\infty)\). To estimate the remaining component in each closed layer, we use \(u_p=(U+J_p-C_{21}u_q)/C_{22}\). The lift adds only \(C\|u\|_{0,\alpha;\mathrm{pw}}\) to the desired bound. By interpolation, this term is bounded by \(\varepsilon\|u\|_{1,\alpha;\mathrm{pw}}+C_\varepsilon\|u\|_2\), and we absorb the first term into the left-hand side.
\end{proof}

\begin{lemma}\label{lem_schauder_estimate}
Under \cref{def_transmission_system}, for every fixed
\({\tau}\ge0\) there is a constant \(C_{{\tau}}\),
depending only on
 \({\alpha},{p_0},{p_I}, {\rho}_\pm,{\sigma_S},{\sigma_I}, {\nu_0},{K},{\tau},\) 
such that
\begin{equation}\label{eq_transmission_schauder}
\|u\|_{2,{\alpha};\mathrm{pw}}
\le C_{{\tau}}\qty(
\|{\mathscr T}_{{\mathbf C},{\tau}}u\|_{{\widehat Y}}
+\|u\|_{L^2\qty({\Omega})})
\qquad\qty(u\in{\widehat X}).
\end{equation}
The estimate applies in particular to even and odd periodic functions;
under \eqref{eq_transmission_reflection_covariance}, the corresponding
operator and data spaces preserve these symmetries.
\end{lemma}
\begin{proof}
Write \(\mathscr T_{C,0}u=(f_+,f_-,g_S,g_I)\). We first estimate \(u\) by \cref{lem_transmission_gradient} with \(G=0\). Differentiating tangentially, we apply the same lemma to \(v=u_q\) with \(G_v=f e_q-C_q\nabla u\), zero scalar source and boundary data \(\partial_qg_S,\partial_qg_I\). The effective normal flux of the differentiated equation is the tangential derivative of the original normal flux. In each layer
\begin{equation}\label{eq_transmission_normal_second}
 C_{22}u_{pp}=f-C_{11}u_{qq}-2C_{12}u_{qp}-(\operatorname{div}C)\cdot\nabla u.
\end{equation}
Division by \(C_{22}\ge\nu_0\) controls the remaining derivative. Adding a fixed mass shift contributes \(\tau\rho u\) to \(f\); interpolation against \(\|u\|_2\) absorbs its layerwise \(C^\alpha\) norm.
\end{proof}

\begin{lemma}\label{lem_shifted_isomorphism}
For the system of \cref{def_transmission_system}, we can choose
\({\tau}_0>0\), depending only on
 \({p_0},{p_I},{\rho}_\pm, \allowbreak{\sigma_S},{\sigma_I},{\nu_0},\) 
so that, for every fixed \({\tau}\ge{\tau}_0\),
 \({\mathscr T}_{{\mathbf C},{\tau}}: {\widehat X}\longrightarrow{\widehat Y}\) 
is a bounded isomorphism and
\begin{equation}\label{eq_transmission_inverse_estimate}
\|{\mathscr T}_{{\mathbf C},{\tau}}^{-1}F
\|_{2,{\alpha};\mathrm{pw}}
\le C_{{\tau}}'\|F\|_{{\widehat Y}}.
\end{equation}
Under \eqref{eq_transmission_reflection_covariance}, the same
isomorphism and inverse estimate hold separately on the even and odd
periodic subspaces.
\end{lemma}
\begin{proof}
On \(V_2\) the weak form and load are
\begin{equation}\label{eq_transmission_weak_form}
 B_\tau(u,v)=\sum_i\int_{\Omega_i}(C_i\nabla u\cdot\nabla v+\tau\rho_iuv)
 -\sigma_S\int u_Sv_S+\sigma_I\int u_Iv_I ,
\end{equation}
and \(\mathcal L(v)=-\sum_i\int f_iv+\int g_Sv_S-\int g_Iv_I\). The trace inequality bounds the sum of squared trace norms by \(\epsilon\|\nabla u\|_2^2+C_\epsilon\|u\|_2^2\). We choose \(\epsilon\) so that the boundary loss is at most \(\nu_0/4\) in the gradient term, then take \(\tau\min\rho_i\) larger than the remaining mass loss by one. The Lax--Milgram theorem now gives a unique weak solution satisfying \(\|u\|_{H^1}\le C\|F\|_{\widehat Y}\).

For smooth data, tangential difference quotients show inductively that \(u_j=\partial_q^j u\in V_2\): in the differentiated form the extra loads are \(-\sum_i\sum_{r=1}^j\binom jr\int(\partial_q^r C_i)\nabla u_{j-r}\cdot\nabla v\), bounded by the preceding estimates. Equation \eqref{eq_transmission_normal_second} with source \(f+\tau\rho u\) now recovers all normal orders by induction, hence smoothness on each closed layer. We may therefore apply the Schauder estimate to this solution.

For general data we extend and mollify the coefficients and data on each strip separately. In doing so we preserve symmetry, ellipticity at least \(\nu_0/2\) and uniform bounds in the stated H\"older norms. Coefficients converge in \(C^1\), sources in \(C^0\), and boundary data in \(C^1\). The coercivity margin and the Schauder bound give a uniform piecewise \(C^{2,\alpha}\) estimate for the approximants. After passage to a subsequence, we obtain convergence in \(C^2\) on each closed layer. The limit satisfies the equations with the same seminorm bound. Uniqueness of the weak solution identifies this limit with the inverse image of the data. If the operator commutes with reflection, uniqueness also shows that the inverse preserves parity.
\end{proof}

\begin{proposition}\label{prop_transmission_fredholm}
Under \cref{def_transmission_system},
 \({\mathscr T}_{{\mathbf C},0}: {\widehat X}\longrightarrow{\widehat Y}\) 
is Fredholm of index zero.
Under \eqref{eq_transmission_reflection_covariance}, its restrictions
to the even and odd periodic subspaces are also Fredholm of index zero.
The allowed Robin data include
 \({\sigma_S}={\rho}_+{g}, \qquad {\sigma_I}={\Delta\rho}{g}.\) 
\end{proposition}
\begin{proof}
Set \(\mathscr Ju=(\rho_+u_+,\rho_-u_-,0,0)\). The layerwise compact embedding \(C^{2,\alpha}\hookrightarrow C^{0,\alpha}\) makes \(\mathscr J:\widehat X\to\widehat Y\) compact. Since \(\mathscr T_{\mathbf C,0}=\mathscr T_{\mathbf C,\tau}+\tau\mathscr J\), \cref{lem_shifted_isomorphism} gives index zero. Under reflection covariance the shifted isomorphism restricts to either parity class, and \(\mathscr J\) preserves both, so the same argument applies there.
\end{proof}

\subsection{The nonlinear operator and controlled compactness}\label{subsec_nonlinear_compactness}

\begin{proposition}\label{prop_height_fredholm}
The domain \(\mathcal U\) is open and convex, and \(\mathcal F:\mathbb R\times\mathcal U\to Y\) is real analytic. Write
\begin{equation}\label{eq_height_coefficient_matrices}
 A_i=\begin{pmatrix} h_p^{-1}&-h_qh_p^{-2}\\-h_qh_p^{-2}&(1+h_q^2)h_p^{-3}\end{pmatrix},\qquad
 \mathbf C_i=\rho_i A_i .
\end{equation}
With \(j_i=A_i\nabla m\), the height and Bernoulli derivatives are
\begin{equation}\label{eq_height_complete_linearization}
 L_hm=(\operatorname{div}j_+,\operatorname{div}j_-,j_{+,p}|_S-gm_S,[\rho j_p]_I-\Delta\rho g m_I),
\end{equation}
Here \(j_p\) denotes the vertical component, and \(\mathcal F_Q=(0,0,1/2,\Delta\rho/2)\).
For every \({h}\in{\mathcal U}\),
\({L}_{{h}}:{X}\to{Y}\)
is Fredholm of index zero.
For every \(R\ge1\), \(\varepsilon>0\), the bounds
\(\|{h}\|_{{X}}\le R\) and
\(\min{h}_p\ge\varepsilon\) imply
\begin{equation}\label{eq_height_linear_estimate}
\|m\|_{{X}}
\le C\qty(
\|{L}_{{h}}m\|_{{Y}}
+\|m\|_{L^2\qty({\Omega})}),
\end{equation}
where \(C\) depends only on \(R,\varepsilon\) and the fixed data.
\end{proposition}
\begin{proof}
The condition \(h_p>0\) defines an open convex set of heights. In each layer the reciprocal of \(h_p\) is analytic by its local Neumann series in the \(C^{1,\alpha}\) algebra. Taking products, derivatives and traces gives the asserted analyticity. Differentiation yields the stated linearization. Since the height is even, both diagonal coefficients are even and the off-diagonal coefficient is odd. Under the stipulated bounds \(\det A_i=h_p^{-4}\ge R^{-4}\), while \(\operatorname{tr}A_i\) and \(\|A_i\|_{1,\alpha}\) are bounded in terms of \(R,\varepsilon\). Multiplying the two interior residuals and the surface residual by \(\rho_+,\rho_-,\rho_+\), respectively, converts \(L_h\) to \(\mathscr T_{\mathbf C,0}\) with Robin constants \(\rho_+g,\Delta\rho g\). This invertible rescaling, \cref{lem_schauder_estimate,prop_transmission_fredholm} and reflection covariance give the assertion on the even spaces.
\end{proof}

\begin{proposition}\label{prop_regularity_compactness}
For every \(R\ge1\), \(\varepsilon>0\), there is a constant depending
only on
 \(R,\varepsilon,{\alpha},{p_0},{p_I}, {g},{\rho}_\pm,{\gamma}_\pm\) 
such that each solution with
 \({\mathcal F}\qty({Q},{h})=0,\qquad |{Q}|+\|{h}\|_{{X}}\le R,\qquad \min{h}_p\ge\varepsilon\) 
satisfies
\begin{equation}\label{eq_height_third_derivative_bound}
\max_{i=\pm}
\|{h}_i\|_{C^{3,{\alpha}}\qty(\overline{{\Omega}_i})}
\le C.
\end{equation}
Moreover,
\begin{equation}\label{eq_height_compact_zero_set}
\begin{aligned}
{\mathscr Z}_{R,\varepsilon}
=\{\qty({Q},{h})\in
\mathbb R\times{\mathcal U}:\
&{\mathcal F}\qty({Q},{h})=0,\\
&|{Q}|+\|{h}\|_{{X}}\le R,\
\min{h}_p\ge\varepsilon\}
\end{aligned}
\end{equation}
is compact in \(\mathbb R\times{X}\).
\end{proposition}
\begin{proof}
Tangential differentiation of the weak equations gives \(L_h\phi=0\) for \(\phi=h_q\in V_2\), odd in \(q\), with common interface trace. A uniform coercive shift exists by \cref{prop_height_fredholm}. The shifted data are \((-\tau\rho_+\phi_+,\allowbreak-\tau\rho_-\phi_-,0,0)\), bounded in the odd H\"older target. Weak uniqueness and the inverse from \cref{lem_shifted_isomorphism} identify the strong solution with \(\phi\), hence \(\|h_q\|_{2,\alpha;\mathrm{pw}}\le C(R,\varepsilon)\). In each layer
\[
 h_{pp}=(2h_ph_qh_{qp}-h_p^2h_{qq}-\gamma_i h_p^3)/(1+h_q^2).
\]
By the estimate for \(h_q\), the normal derivative of the right-hand side is bounded in \(C^\alpha\). We thus obtain the required bound for \(h_{ppp}\) as well. The resulting bounded set in piecewise \(C^{3,\alpha}\) is precompact in \(X\); the margin on \(h_p\), the traces, the norm constraints and the zero equation pass to the limit.
\end{proof}

\section{Local and global bifurcation}\label{sec_bifurcation}
We now combine the spectral calculation with the Fredholm property to obtain local bifurcation. We then continue the resulting branch in the physical solution space.

\subsection{The analytic local branch}\label{subsec_local_bifurcation}

\begin{proof}[Proof of \cref{thm_local_bifurcation}]
The change \((Q,h)\leftrightarrow(\kappa,h)\) in \eqref{eq_gauge} is an affine analytic isomorphism on the entire physical open set. We set \(G(\kappa,h)=\mathcal F(\kappa+2g\overline h_S,h)\). Its state derivative is
 \(D_hG\,m=L_hm+2g\overline m_S(0,0,1/2,\Delta\rho/2).\) 
The rank-one perturbation preserves Fredholm index zero by \cref{prop_height_fredholm}. At a laminar solution the mean fluxes \(j_i=a_i^3m_p\) of a mean kernel element are constant. Its surface equation gives \(j_+=0\), whence \(m_S=m_I\). The interface equation becomes
\(\rho_+j_+-\rho_-j_-+\Delta\rho g(m_S-m_I)=0\); therefore \(j_-=0\) and the bed trace gives \(m=0\). This holds even at a turning point of \(Q(\lambda)\). At fixed \(Q\), \cref{lem_zero_mode} still applies.

Every nonzero Fourier mode has zero surface mean. At a positive fundamental root we have \(\nu_1=-1\). Hence \(-k^2<\nu_1\) excludes every \(k\ge2\). The kernel of the translated problem \(G(\kappa,H(\kappa)+w)=0\) is thus spanned by \(\varphi_*\). The Green functional in \eqref{eq_green_range_identity} annihilates the additional constant term. It spans the cokernel, and its pairing with the derivative of the translated linearization applied to \(\varphi_*\) is \(\pi\rho_+\mathscr D_1'/\mathsf D_1\) by \eqref{eq_exact_transversality_pairing}. In particular it is nonzero.

For completeness, write \(w=s\varphi_*+v\) with \(\mathfrak b(v)=0\). Projection onto the range and the analytic implicit function theorem determine \(v=v(\kappa,s)\), with \(v(\kappa,0)=0\) and \(v_s(\lambda_*,0)=0\). The remaining scalar equation vanishes at \(s=0\). After dividing analytically by \(s\), we obtain an equation whose \(\kappa\) derivative at the base point is the nonzero Green pairing above (the range correction is annihilated). A second application of the implicit function theorem gives exactly one analytic curve. Translation by \(\pi\) preserves the residual and changes the sign of the amplitude, so uniqueness gives the parity assertions and the stated expansions. Reconstructing \(Q\) by \eqref{eq_gauge} proves local completeness in the physical coordinates.
\end{proof}

\begin{remark}\label{prop_nodal_scope}
At fixed \(p>p_0\), \(h_q(s)(\pi/2,p)=-sM_*(p)+O(s^2)\) distinguishes the two local orientations. Even periodic symmetry gives \(h_q=h_{qp}=0\) on the axes, including both one-sided interface traces.
\end{remark}

\subsection{The physical global continuum}\label{subsec_global_bifurcation}

For related compact and analytic continuation theorems see \cite[Theorem~1.3]{Rabinowitz1971} and \cite[Chapter~9]{BuffoniToland2003}. We use the Fredholm component theorem of Shi and Wang \cite[Theorem~4.3]{ShiWang2009}.

\begin{theorem}\label{thm_global_bifurcation}
Assume the data and root conditions of \cref{thm_local_bifurcation}.
The base point \(z_*\) belongs to \(\mathscr S\). Let \(\mathcal C\) be its component in
\cref{def_nodal_set}.
It contains the local curve and satisfies at least one alternative: (i) \(|Q|+\|h\|_X\) is unbounded, (ii) it contains a laminar point with \(\lambda_1\ne\lambda_*\), or (iii) it contains a sequence with \(\min h_p\to0\). In (ii) the Bernoulli values may coincide.
The third alternative concerns the boundary of \(\mathcal U\). On a height-norm ball of radius \(R\) the surface law gives \(Q-2gh(q,p)\ge Q-2gh(q,0)\ge R^{-2}\).
\end{theorem}
\begin{proof}
The local curve lies in \(\mathscr S\), and its base component is \(\mathcal C\). If (i) and (iii) both fail, \cref{prop_regularity_compactness} places this relatively closed component in a compact subset of the physical open set. Suppose also that its only laminar flow is \(z_*\). Work in the global coordinates \((\kappa,h)\) from \eqref{eq_gauge}. We choose open intervals \(J\Subset J'\Subset\mathcal I\) containing \(\lambda_*\) and a smooth cutoff equal to one on \(\overline J\) and supported in \(J'\). With this cutoff take the convex combination of \(H(\kappa)\) and the bed-zero reference \(p-p_0\), extending by the latter outside \(J'\); denote the resulting horizontal profile by \(\widetilde H(\kappa)\). Its vertical derivative stays positive.

Delete from the physical open set the laminar zeros with \(\kappa\notin J\), and the points with \(\kappa\notin J,\ h=\widetilde H(\kappa)\). These sets are relatively closed: horizontality and the zero equation are closed conditions, and the reference graph is closed. Neither meets \(\mathcal C\), since a horizontal solution is laminar. The connected component of the remaining open set containing \(\mathcal C\) is open (Banach spaces are locally path connected). There change variables by \(w=h-\widetilde H(\kappa)\). Every point of the trivial line belonging to this open domain is a zero, and all its laminar zeros lie on that line.

Throughout this open domain, at fixed \(\kappa\), the derivative in \(w\) sends \(m\) to \(L_hm+2g\overline m_S\partial_Q\mathcal F\), with \(h=w+\widetilde H(\kappa)\). By \cref{prop_height_fredholm} its first term is Fredholm of index zero and its second has rank at most one. The intersection of the trivial line with the domain consists exactly of \((\kappa,0)\) with \(\kappa\in J\): the graph over \(J\) is connected and contains the base, and its complement on the line was deleted. This verifies the index condition throughout the domain. At the base flow, the kernel is simple and the crossing is nonzero by the local calculation. The residual is smooth and its zero set is locally compact by \cref{prop_regularity_compactness}. In this domain the closure of the nontrivial zeros is the restriction of \(\mathscr S\), since taking closure commutes with restriction to an open set. Its base component is still \(\mathcal C\): it contains that entire connected set and any larger connected set would contradict maximality in \(\mathscr S\). The component alternative of Shi and Wang \cite[Theorem~4.3]{ShiWang2009} now precludes a compact base component with no other trivial point. This gives the asserted alternatives. To obtain the pressure-head bound, we use the surface equation and \(h_p>0\).
\end{proof}

Here a return means that the connected set contains a second physical laminar point. This formulation requires no global parametrization.
\section{Nodal geometry and the monotone continuum}\label{sec_nodal_geometry}
In this section we follow the sign of the streamline slope from the local branch along the global component. The main point is to control this sign where a symmetry axis meets the density interface.

\subsection{Nodal orientation and strict derivative bounds}\label{subsec_nodal_positivity}

\begin{theorem}\label{thm_nodal_positivity}
For every \(\qty({Q},{h})\in{\mathscr P}\) there is
\(\chi_{{h}}>0\) such that
\begin{equation}\label{eq_nodal_weighted_positivity}
-{h}_q\qty(q,p)\ge
\chi_{{h}}\sin q\,\qty(p-{p_0})
\qquad\qty(0\le q\le\pi,\ {p_0}\le p\le0).
\end{equation}
Every interior streamline, the free surface and the interface are
strictly decreasing for \(0<q<\pi\). All layerwise one-sided traces satisfy
\begin{equation}\label{eq_nodal_boundary_derivatives}
\begin{aligned}
&{h}_{qq}\qty(0,p)<0,\qquad
{h}_{qq}\qty(\pi,p)>0
&&{p_0}<p\le0,\\
&{h}_{qp}\qty(q,{p_0})<0
&&0<q<\pi,\\
&{h}_{qqp}\qty(0,{p_0})<0,\qquad
{h}_{qqp}\qty(\pi,{p_0})>0,\\
&{h}_{qp}\qty(0,p)={h}_{qp}\qty(\pi,p)=0
&&p\in\overline{{I}_\pm}.
\end{aligned}
\end{equation}
The set \({\mathscr P}\) is both relatively open and relatively closed
in \({\mathscr Z}\setminus{\Lambda}\), for the topology of
\(\mathbb R\times{X}\). For openness, convergence of actual
solutions in this topology improves to layerwise \(C^3\) by
\cref{prop_regularity_compactness}. Every connected component of the
nonlaminar solution set meeting \({\mathscr P}\) is therefore contained
in \({\mathscr P}\).

Under the root hypotheses and amplitude normalization of
\cref{thm_local_bifurcation}, all sufficiently small positive-amplitude
members of its local curve belong to \({\mathscr P}\). Half-period
translation gives the opposite nodal direction for negative amplitude.
\end{theorem}
\begin{proof}
We adapt the maximum-principle argument of Constantin and Strauss
\cite{ConstantinStrauss2004} to the transmission equation for the
negative tangential derivative \({W}\) from
\cref{def_nodal_set}, using local barriers where the vertical axes
meet the interface, surface and bed. From the signs of the derivatives we obtain a weighted lower bound. We then use uniform convergence to prove relative closedness. For relative openness and the sign on the positive local branch, we integrate from the axes and the bed and apply the weighted bound.

\smallskip\noindent\textit{Step 1. Positivity and the ordinary Hopf signs.}\enspace
We prove positivity in both layers and on the interface and surface
for \(0<q<\pi\), and obtain the ordinary inward derivative signs
away from boundary intersections.
By \cref{prop_regularity_compactness}, each solution is
layerwise \(C^{3,{\alpha}}\). Set
\begin{equation}\label{eq_nodal_tangential_system}
\begin{aligned}
{W}&=-{h}_q,\\
{\mathbf C}_i&={\rho}_i
\begin{pmatrix}
{h}_p^{-1}&-{h}_q{h}_p^{-2}\\
-{h}_q{h}_p^{-2}&
\qty(1+{h}_q^2){h}_p^{-3}
\end{pmatrix},\\
{\beta_S}&={\rho}_+{g},\qquad
{\beta_I}={\Delta\rho}{g}.
\end{aligned}
\end{equation}
The matrices are layerwise \(C^{2,{\alpha}}\), symmetric and
positive definite, with a uniform ellipticity constant on the closed
layers of this fixed solution. Exact tangential differentiation gives
\begin{equation*}
\begin{aligned}
\operatorname{div}\qty({\mathbf C}_i\nabla{W})&=0,\\
\qty({\mathbf C}_+\nabla{W})_p
&={\beta_S}{W}&&p=0,\\
\qty[\qty({\mathbf C}\nabla{W})_p]
&={\beta_I}{W}&&p={p_I},\\
{W}&=0&&q=0,\pi\text{ or }p={p_0},\\
{W}^+&={W}^-&&p={p_I}.
\end{aligned}
\end{equation*}
The function \({W}\) is nonnegative and not identically
zero. The strong maximum principle makes it either positive or
identically zero in each open layer. If one layer were zero and the
other nonzero, Hopf's lemma at any interface point with \(0<q<\pi\)
would give a nonzero conormal derivative, contradicting the zero trace
and the flux jump. Both open layers are therefore strictly positive.

If \({W}\qty(q_0,{p_I})=0\) for
\(0<q_0<\pi\), its common tangential trace has a minimum there, so
\(\partial_q{W}=0\). Hopf's lemma on both sides gives
 \(\partial_p{W}\qty(q_0,{p_I}+)>0, \qquad \partial_p{W}\qty(q_0,{p_I}-)<0.\) 
The flux jump is then strictly positive, whereas its prescribed value
\({\beta_I}{W}\) is zero. At a zero on the
free surface, the tangential derivative again vanishes and Hopf's
lemma gives \(\partial_p{W}<0\), contradicting zero
surface flux. Thus \({W}>0\) for
\(0<q<\pi\), \({p_0}<p\le0\), including the interface and
surface. Ordinary Hopf comparison at the bed gives
\(\partial_p{W}\qty(q,{p_0})>0\).
On the vertical edges away from their intersections with the bed,
interface and surface, it gives
 \(\partial_q{W}\qty(0,p)>0,\qquad \partial_q{W}\qty(\pi,p)<0.\) 

\smallskip\noindent\textit{Step 2. Comparison on small rectangles.}\enspace
For the comparisons at the interface and surface, we absorb the trace
terms into the elliptic energy. This gives a comparison principle for
barriers satisfying the conormal inequalities and the Dirichlet
comparison on the edges where values are prescribed.
Consider a rectangle of width \(\epsilon\), with height
at most \(2\epsilon\) in each layer, possibly meeting the interface or
the free surface. The vertical symmetry edge and the artificial edges
carry Dirichlet comparison data. Let \(v\) be the negative part of a
comparison difference, vanishing at \(x=0\) and on the artificial
Dirichlet edges. Horizontal Poincar\'e and one-dimensional trace
estimates give
\begin{equation*}
\begin{aligned}
\int v^2&\le\epsilon^2\int|v_x|^2,\\
\int_{\text{horizontal natural edges or interface}}v^2
&\le C\qty(\epsilon^{-1}\int v^2+\epsilon\int|v_p|^2)
\le C\epsilon\int|\nabla v|^2.
\end{aligned}
\end{equation*}
The shapes below are \(\epsilon\times\epsilon\) or
\(\epsilon\times2\epsilon\), so \(C\) is a fixed shape constant.
Let \({\nu_0}>0\) be a lower ellipticity bound for the fixed
matrices. We choose \(\epsilon\) small enough that
 \(C\epsilon\qty({\beta_S}+|{\beta_I}|) <{\nu_0}/2.\) 
The local weak form is coercive:
\begin{equation*}
\sum_i\int{\mathbf C}_i\nabla v\cdot\nabla v
-{\beta_S}\int_Sv^2
+{\beta_I}\int_Iv^2
\ge\frac{{\nu_0}}2\int|\nabla v|^2.
\end{equation*}
A natural boundary absent from the rectangle has an empty integration
region. Consequently, a barrier \(\zeta\) with
\begin{equation*}
\begin{aligned}
\operatorname{div}\qty({\mathbf C}_i\nabla\zeta)&\ge0,\\
\qty({\mathbf C}_+\nabla\zeta)_p
-{\beta_S}\zeta&\le0&&\text{if the surface is present},\\
\qty[\qty({\mathbf C}\nabla\zeta)_p]
-{\beta_I}\zeta&\ge0&&\text{if the interface is present}
\end{aligned}
\end{equation*}
satisfies \({W}\ge\varepsilon\zeta\) whenever this
inequality holds on the Dirichlet edges. Indeed, for nonnegative \(v\) and \(J=\mathbf C\nabla\zeta\), the boundary terms of the barrier form are \(+\int_S(J_p-\beta_S\zeta)v\) and \(-\int_I([J_p]-\beta_I\zeta)v\). Both terms are nonpositive. The same is true of the interior term \(-\int \operatorname{div}(J)v\). Thus \(W-\varepsilon\zeta\) is a weak supersolution. We obtain the comparison by testing with the negative part and applying coercivity.

\smallskip\noindent\textit{Step 3. The side derivatives at the interface and surface.}\enspace
We use Step 2 to prove strict positivity of the inward derivative
where either vertical side meets the interface or the surface.
Near \(q=0\) put \(x=q\); near \(q=\pi\) put
\(x=\pi-q\) and reflect the matrix accordingly. Denote the transformed
matrix by the same symbol and set
 \(d_1=\partial_x\qty({\mathbf C})_{11} +\partial_p\qty({\mathbf C})_{12},\qquad d_2=\partial_x\qty({\mathbf C})_{12} +\partial_p\qty({\mathbf C})_{22}.\) 
Evenness makes the diagonal coefficients even in \(x\) and the
off-diagonal coefficient odd. On both closed layer rectangles, a
finite constant \(C_0\), controlled by the fixed solution's
\(C^{3,{\alpha}}\) regularity, therefore satisfies
 \(|\qty({\mathbf C})_{12}|\le C_0x,\qquad |d_1|\le C_0x,\qquad |d_2|\le C_0,\qquad \qty({\mathbf C})_{22}\ge{\nu_0}>0.\) 
In \(0<x<\epsilon\), \(|p-{p_I}|<\epsilon\), choose
\(\epsilon\) small enough for the comparison principle and set
 \(\zeta=x\exp\qty(b|p-{p_I}|),\qquad \sigma=\operatorname{sign}\qty(p-{p_I}).\) 
In each layer the complete expansion gives
\begin{equation*}
\begin{aligned}
\frac{\operatorname{div}\qty({\mathbf C}_i\nabla\zeta)}{\zeta}
&=\qty({\mathbf C}_i)_{22}b^2
+\sigma b\qty(2\frac{\qty({\mathbf C}_i)_{12}}x+d_2)
+\frac{d_1}x\\
&\ge{\nu_0}b^2-3C_0b-C_0.
\end{aligned}
\end{equation*}
At the interface \(\zeta=x\), and the conormal remainder is
\begin{equation*}
\begin{aligned}
\frac{\qty[\qty({\mathbf C}\nabla\zeta)_p]
-{\beta_I}\zeta}{\zeta}
&=b\qty(\qty({\mathbf C}^+)_{22}
+\qty({\mathbf C}^-)_{22})
+\frac{\qty({\mathbf C}^+)_{12}-\qty({\mathbf C}^-)_{12}}x
-{\beta_I}\\
&\ge2{\nu_0}b-2C_0-{\beta_I}.
\end{aligned}
\end{equation*}
We choose \(b\) large enough to make both lower bounds nonnegative.
On \(x=\epsilon\), \({W}\) is strictly positive,
including its common interface trace. On
\(p={p_I}\pm\epsilon\), the quotient
\({W}/x\) extends to a strictly positive value by
ordinary side Hopf comparison. Hence a sufficiently small
\(\varepsilon>0\) gives the Dirichlet comparison on all artificial
edges. The comparison principle yields
 \({W}\ge\varepsilon\zeta,\qquad \partial_x{W}\qty(0,{p_I}) \ge\varepsilon>0.\) 
This is the common tangential derivative of the two interface traces.

For the surface intersection, use
\(0<x<\epsilon\), \(-\epsilon<p<0\) and
\(\zeta=x\exp\qty(-b p)\). The full interior and boundary estimates are
\begin{equation*}
\begin{aligned}
\frac{\operatorname{div}\qty({\mathbf C}_+\nabla\zeta)}{\zeta}
&=\qty({\mathbf C}_+)_{22}b^2
-b\qty(2\frac{\qty({\mathbf C}_+)_{12}}x+d_2)
+\frac{d_1}x\\
&\ge{\nu_0}b^2-3C_0b-C_0,\\
\frac{\qty({\mathbf C}_+\nabla\zeta)_p
-{\beta_S}\zeta}{\zeta}
&=\frac{\qty({\mathbf C}_+)_{12}}x
-b\qty({\mathbf C}_+)_{22}-{\beta_S}\\
&\le C_0-{\nu_0}b-{\beta_S}.
\end{aligned}
\end{equation*}
For large \(b\), the first is nonnegative and the second
nonpositive. Strict positivity at the surface controls
\(x=\epsilon\); ordinary side Hopf comparison controls
\(p=-\epsilon\). Choosing \(\varepsilon>0\) and applying the comparison principle gives
\(\partial_x{W}\qty(0,0)>0\).
The strict side derivative now holds for all
\({p_0}<p\le0\); reflection gives its negative sign at
\(q=\pi\).

\smallskip\noindent\textit{Step 4. The mixed derivative at the bed corners.}\enspace
We set \(t=p-{p_0}\). At each bed corner we prove that
\(\partial_{xt}{W}\) is positive by Dirichlet comparison
with \(\zeta=xt\exp\qty(b t)\) and integration from both axes.
Use the lower-layer rectangle \(0<x<\epsilon\), \(0<t<\epsilon\)
and the coefficient bounds of Step 3. The complete
expansion is
\begin{equation*}
\begin{aligned}
\frac{\operatorname{div}\qty({\mathbf C}_-\nabla\zeta)}{x\exp\qty(b t)}
={}&\qty({\mathbf C}_-)_{22}\qty(2b+b^2t)
+2\frac{\qty({\mathbf C}_-)_{12}}x\qty(1+b t)\\
&+\frac{d_1}x\,t+d_2\qty(1+b t)\\
\ge{}&{\nu_0}\qty(2b+b^2t)
-3C_0\qty(1+b t)-C_0t.
\end{aligned}
\end{equation*}
We choose \(b\) large enough to make the last expression nonnegative
for \(0\le t\le\epsilon\). On \(x=\epsilon\), ordinary bed Hopf
comparison makes \({W}/t\) strictly positive; on
\(t=\epsilon\), ordinary side Hopf comparison makes
\({W}/x\) strictly positive. A small
\(\varepsilon>0\) gives comparison on these two edges, while both
functions vanish on the other two. Dirichlet comparison yields
\({W}\ge\varepsilon xt\exp\qty(b t)\).
Since \({W}\) is \(C^2\) and vanishes on both axes,
\begin{equation*}
{W}\qty(x,t)
=\int_0^x\int_0^t
\partial_{xt}{W}\qty(u,v)\,dv\,du,
\qquad
\partial_{xt}{W}\qty(0,0)\ge\varepsilon>0.
\end{equation*}
Returning to the original coordinates gives
 \(\partial_{qp}{W}\qty(0,{p_0})>0,\qquad \partial_{qp}{W}\qty(\pi,{p_0})<0.\) 

\smallskip\noindent\textit{Step 5. The weighted lower bound and boundary signs.}\enspace
It remains to make the positivity estimate uniform. In the open half-strip the quotient \(W/(\sin q\,(p-p_0))\) is positive. At the vertical edges and bed it has positive one-sided limits by the respective inward derivatives; at the bed corners use the double integral of \(W_{qp}\) and its strict corner sign. Step 3 gives the corresponding bounds at the interface and surface corners. By compactness we obtain a positive minimum \(\chi_h\). Evenness gives \(h_q=h_{qp}=0\) on the axes.

\smallskip\noindent\textit{Stability and cone entry.}
Uniform convergence of \(h_q\) gives relative closedness among nonlaminar waves. For solutions converging in \(X\), \cref{prop_regularity_compactness} and uniqueness of the limit give layerwise \(C^3\) convergence. If \(u\) vanishes on the axes and bed and has a common tangential derivative at the interface, integration from the nearest axis and from the bed gives
\[
 |u(q,p)|\le(\pi/2)\sin q\,(p-p_0)\|u_{qp}\|_{\infty;\mathrm{pw}} .
\]
Applying this estimate to the difference of the slopes, we see that the weighted error is eventually less than \(\chi_h/2\), proving openness. On the local curve \(W(s)/s\to M_*\sin q\) in layerwise \(C^{1,\alpha}\). Applying the odd shifted inverse to this quotient gives a uniform layerwise \(C^{2,\alpha}\) bound, since its shifted interior data are \(-\tau\rho_iW(s)/s\), its boundary data vanish and the coefficient bounds are uniform. Thus convergence holds in layerwise \(C^2\). Positivity of \(M_*\) and ODE uniqueness at its bed zero give \(\inf M_*(p)/(p-p_0)>0\). The same double integration now shows entry into the positive cone for \(s>0\); half-period translation gives the other orientation.
\end{proof}

\subsection{The global continuum with one crest and one trough}\label{subsec_monotone_continuum}

\begin{theorem}\label{thm_unilateral_global}
Under the hypotheses of \cref{thm_local_bifurcation}, the small positive half-branch belongs to \(\mathcal C^+\); every nonlaminar member satisfies \cref{thm_nodal_positivity}. The three alternatives of \cref{thm_global_bifurcation} hold with \(\mathcal C^+\) in place of \(\mathcal C\). Half-period translation gives the opposite orientation, and the two components can intersect only at laminar flows.
\end{theorem}
\begin{proof}
Write \(P=\mathscr S^+\), \(T=T_\pi\). The sign on the local branch and relative closedness of the nodal condition give the first assertions and \(P\setminus\Lambda=\mathscr P\). Evenness gives \((Th)_q(q,p)=-h_q(\pi-q,p)\) on the half-period, so \(P\cap TP\subset\Lambda\). Nodal openness supplies a neighborhood of each point of \(P\setminus\Lambda\) in \(\mathscr S\) lying wholly in \(P\); the reverse orientation has the same property. The local description of all solutions and the sign of the bifurcating slope at \(z_*\) put a neighborhood in \(P\cup TP\) and isolate \(z_*\) among the laminar points of \(\mathscr S\).

Suppose, to the contrary, that none of the three alternatives holds. Then \(\mathcal C^+\) is compact by \cref{prop_regularity_compactness} and contains only one laminar point. We choose a bounded open neighborhood \(O\) of the component whose closure has \(\min h_p\ge\varepsilon>0\) and avoids \((P\cap\Lambda)\setminus\{z_*\}\), a closed set by the local isolation. In the compact space \(K=P\cap\overline O\), the component \(\mathcal C^+\) misses \(K\cap\partial O\). In a compact Hausdorff space each component is a quasicomponent. We may therefore separate each boundary point by a clopen set containing the component, then intersect finitely many by compactness of the boundary. The resulting set \(K_1\subset O\) is compact, relatively open in \(P\), and meets \(\Lambda\) only at \(z_*\). Now \(K_1\cup TK_1\) is compact and relatively open in \(\mathscr S\): away from \(z_*\) use nodal openness, and at \(z_*\) use the two local halves and the relative openness of each piece. It contains the bilateral component and no other laminar point. This contradicts \cref{thm_global_bifurcation}.
\end{proof}

The closure defining each component can contain laminar points. Normalizing the tangential derivative along a nonlaminar sequence identifies the admissible dispersion modes at those points.

\begin{proposition}\label{prop_laminar_return}
Consider convergence in \(\mathbb R\times X\) to an admissible laminar
solution \(\qty(Q\qty(\lambda),H\qty(\lambda))\). If the approximating solutions
are nonlaminar, then \(\mathscr D_k\qty(\lambda)=0\) for some \(k\ge1\).
If they also satisfy \(h_q\le0\) for \(0<q<\pi\) in both layers,
then
\begin{equation}\label{eq_nodal_return_fundamental}
\nu_1\qty(\lambda)=-1,\qquad
\mathscr D_1\qty(\lambda)=0,\qquad \mathsf A_1\qty(\lambda)>0.
\end{equation}
\end{proposition}
\begin{proof}
We set \(W=-h_q\). Tangential differentiation of the weak height equations, with the common interface trace and zero bed trace, gives for every \(\phi\in V_2\)
\begin{equation}\label{eq_laminar_return_weak_slope}
 B_h(W,\phi)=\sum_i\int_{\Omega_i}\mathbf C_i(h)\nabla W\cdot\nabla\phi
 -g\rho_+\int W_S\phi_S+g\Delta\rho\int W_I\phi_I=0.
\end{equation}
We normalize \(v_n=W_n/\|W_n\|_2\). The matrices converge uniformly to \(\mathbf C_i(H)\) and share an ellipticity constant. Testing with \(v_n\), the trace inequality \(\|v_n(\cdot,t)\|_2^2\le\epsilon\|\nabla v_n\|_2^2+C_\epsilon\|v_n\|_2^2\) at \(t=0,p_I\) yields a uniform \(H^1\) bound after absorption. A subsequence converges weakly in \(H^1\), strongly in \(L^2\), to an odd unit vector \(v\). Weak convergence of the traces gives \(B_H(v,\cdot)=0\). Testing with \(N(p)\sin kq\), \(N\in V_1\), shows that every nonzero sine coefficient is a vertical eigenfunction of eigenvalue \(-k^2\); at least one is nonzero. In the oriented case \(v\ge0\) on the half-period and its first sine coefficient is nonnegative and nonzero. Orthogonality to the strictly positive principal eigenfunction excludes every higher vertical eigenspace. Now apply \cref{prop_mode_matrix}.
\end{proof}
\begin{corollary}\label{cor_fold}
If the whole admissible set \(\mathcal I\) has just one positive fundamental root and this root is transverse, then
 \(\mathcal C=\mathcal C^+\cup T_\pi\mathcal C^+,\qquad \mathcal C^+\cap T_\pi\mathcal C^+=\{z_*\}.\) 
In particular no further laminar flow, even one with a higher-mode resonance, lies on \(\mathcal C\). Each half has a sequence as in \eqref{eq_oriented_horizontal_stagnation}.
\end{corollary}
\begin{proof}
We set \(P=\mathscr S^+\), \(T=T_\pi\), \(A=P\cup TP\). By \cref{prop_laminar_return} both closed halves have only \(z_*\) as a laminar point. The opposite signs of \(h_q\) exclude a common nonlaminar point. Nodal openness gives a relative neighborhood in \(\mathscr S\) within the corresponding half at every nonlaminar point of \(A\). At \(z_*\) local completeness and cone entry give a whole neighborhood in \(A\). Thus \(A\) is clopen in \(\mathscr S\), so \(\mathcal C\subset A\). The map equal to the identity on \(P\) and to \(T\) on \(TP\) is continuous on \(A\) by closed pasting, since it fixes their common point. Its image of \(\mathcal C\) is connected in \(P\) and contains \(z_*\); hence that image lies in \(\mathcal C^+\). As \(T^2=1\), this proves the forward inclusion. Both \(\mathcal C^+\) and \(T\mathcal C^+\) are connected in \(\mathscr S\) and contain \(z_*\), proving the reverse inclusion. We apply \cref{thm_oriented_stagnation} and translation invariance of \(M\).
\end{proof}

\section{Physical alternatives along the global components}\label{sec_physical_alternatives}

We first obtain a lower bound for the distance between the layers' outer boundaries and the interface, and upper bounds for both speeds, assuming that \(M\) and \(Q\) are bounded. A vertical trace inequality and the large-\(Q\) thinning estimate then control \(Q\) in terms of \(M\) on the full closure of nonlaminar solutions.

\subsection{Circulation and compactness}\label{subsec_physical_compactness}

Bounds on \(M\) and the upper-side interface speed \(V_I\) control
\(Q\) by the free-surface pressure law and the upper-layer circulation.
Pressure continuity transfers the speed bound to the lower layer.
The bounds in both layers then give compactness through closed-layer
regularity.

\begin{proposition}\label{prop_physical_speed_bounds}
Use the quantities of \cref{def_physical_controls}.
Then
\(Q>0\), \(H_{\max}\le |p_0|M\), and
\begin{equation}\label{eq_physical_circulation}
\int_S|\nabla\psi|\,ds
=\int_I|\nabla\psi_+|\,ds+\gamma_+|D_+|.
\end{equation}
In particular,
\begin{equation}\label{eq_physical_bernoulli_bound}
\sqrt{\qty(Q-2gH_{\max})_+}
\le V_I\qty(1+MV_I)+|\gamma_+|\ell_+M .
\end{equation}
The speed bounds on the closed layers are
\begin{equation}\label{eq_physical_layer_speed_bounds}
\begin{aligned}
\|\nabla\psi\|_{L^\infty\qty(D_+)}^2&\le\max\{Q,V_I^2\},\\
\|\nabla\psi\|_{L^\infty\qty(D_-)}&\le
\max_I|\nabla\psi_-|+2|\gamma_-|H_{\max}.
\end{aligned}
\end{equation}
On the interface, pressure continuity gives
\begin{equation}\label{eq_physical_interface_pressure}
|\nabla\psi_-|^2=r|\nabla\psi_+|^2+
\delta\qty(Q+2\gamma_+p_I-2gh\qty(q,p_I)).
\end{equation}
\end{proposition}
\begin{proof}
The surface Bernoulli equation turns the upper-layer circulation into a
bound on \(Q\). We then apply the maximum principle to bound the speed throughout the upper layer.
Continuity of pressure yields the corresponding estimate on the lower side of the interface,
and reflection across the bed gives the lower-layer bound.

The chain rule applied to \(p=-\psi\) gives
\(\Delta\psi=-\gamma_i\) with closed-layer \(C^{2,\alpha}\) regularity and
traces \(0,\ell_+,\ell_++\ell_-\) on \(S,I,B\), respectively.
The kinematic and surface identities are
\begin{equation}\label{eq_physical_kinematics}
\begin{aligned}
\psi_y\circ\Phi_h&=-h_p^{-1},&
|\nabla\psi|^2\circ\Phi_h&=\frac{1+h_q^2}{h_p^2},\\
|\nabla\psi|^2\circ\Phi_h\qty(q,0)&=Q-2gh\qty(q,0).
\end{aligned}
\end{equation}
Integration from the bed gives \(0<h\qty(q,0)\le |p_0|M\), so the
surface identity implies \(Q>0\).

The wave-frame velocity \(\qty(\psi_y,-\psi_x)\) has curl \(\gamma_+\)
and points toward decreasing \(q\) along both upper-layer boundaries.
The positively oriented surface runs from right to left, and the
interface from left to right. Stokes' theorem over one period gives
\eqref{eq_physical_circulation}, with the vertical sides canceling by
periodicity. The interface slope satisfies \(|h_q\qty(q,p_I)|\le MV_I\).
Thus
 \(|S|\ge2\pi,\qquad |I|\le2\pi\qty(1+MV_I),\qquad |D_+|\le2\pi\ell_+M .\) 
The surface speed is at least \(\sqrt{\qty(Q-2gH_{\max})_+}\).
Substituting these estimates into \eqref{eq_physical_circulation} gives
\eqref{eq_physical_bernoulli_bound}.

In the upper layer, \(\Delta|\nabla\psi|^2=2|D^2\psi|^2\ge0\).
The maximum principle with the surface identity gives the upper
bound in \eqref{eq_physical_layer_speed_bounds}. At the interface the
reconstructed pressure is
 \(P_i=\rho_i\qty(Q/2+\gamma_+p_I-gh\qty(q,p_I)-|\nabla\psi_i|^2/2).\) 
Equating the pressures and dividing by \(\rho_-\) gives
\eqref{eq_physical_interface_pressure}.

For the lower-layer speed, the harmonic function
\(w=\psi-\psi|_B+\gamma_-\qty(y+D_0)^2/2\) vanishes at the flat bed.
Its odd reflection is harmonic across the bed, and the reflected
gradient has the same norm as the original gradient at the mirrored
point. We apply the maximum principle to \(|\nabla w|^2\) on this doubled
periodic layer, whose boundary consists of the interface and its
reflection. Then
\begin{equation*}
\begin{aligned}
\max_{\overline D_-}|\nabla\psi|
&\le\max_I|\nabla w|+|\gamma_-|H_{\max}\\
&\le\max_I|\nabla\psi_-|+2|\gamma_-|H_{\max}.
\end{aligned}
\end{equation*}
\end{proof}

\begin{lemma}\label{lem_gradient_to_closed_layer}
Fix \(K\ge1\) and \(0<\alpha<1\). There are \(0<\beta<\alpha\) and
\(C<\infty\), depending only on \(K,\alpha\) and the fixed data, such
that every actual solution with
\begin{equation}\label{eq_height_lipschitz_control}
\max_{i=\pm}\qty(\|h_q\|_{L^\infty\qty(\overline\Omega_i)}
+\|h_p\|_{L^\infty\qty(\overline\Omega_i)}
+\|h_p^{-1}\|_{L^\infty\qty(\overline\Omega_i)})\le K
\end{equation}
satisfies
\begin{equation}\label{eq_gradient_closed_layer_bound}
0<Q,\qquad Q+\max_{i=\pm}\|h\|_{C^{3,\beta}\qty(\overline\Omega_i)}\le C.
\end{equation}
\end{lemma}
\begin{proof}
The matrix in the differentiated equation \eqref{eq_laminar_return_weak_slope} for \(v=h_q\) has determinant \(\rho_i^2h_p^{-4}\). Thus \eqref{eq_height_lipschitz_control} controls both its norm and ellipticity. We apply the lift \eqref{eq_transmission_flux_lifting} to the Robin data: both the lift and its layerwise vertical derivative are bounded by \(C(K)\). The same surface and bed reflections give an interior divergence equation with bounded matrix and bounded sources. On the fixed ball cover the divergence datum lies in \(L^4\), the scalar source in \(L^2\); planar De Giorgi--Nash--Moser estimates \cite[Theorem~8.22]{GilbargTrudinger2001} bound \(\|h_q\|_{C^\beta}\) for some \(0<\beta<\alpha\).

The boundary graphs have bounded \(C^{1,\beta}\) norms. The vertical gaps are bounded below by \(\ell_i/K\), so their Euclidean separation is at least \(\ell_i/(K\sqrt{1+K^2})\) by the Lipschitz bound. Flattening in charts of a fixed smaller radius, the constant Dirichlet data, bounded range of \(\psi\), and \(\Delta\psi=-\gamma_i\) give uniform closed-layer \(C^\beta\) gradient bounds \cite[Theorems~8.32--8.33 and Corollary~8.36]{GilbargTrudinger2001}. Composition with the Lipschitz height map and inversion of \(-\psi_y\ge1/K\) bound \(h_p\) in piecewise \(C^\beta\). We apply \cref{lem_transmission_gradient} at exponent \(\beta\) to \(v\); the coefficients are now piecewise \(C^\beta\). This bounds \(h_{qq},h_{pq}\), and the interior height equation recovers \(h_{pp}\) with the same exponent. The surface law bounds \(Q\le(1+K^2)K^2+2gK|p_0|\); \cref{prop_regularity_compactness} at exponent \(\beta\) gives the third derivatives.
\end{proof}

\begin{proposition}\label{prop_physical_compactness}
For each \(K\ge1\) the solutions with \(M+V_I\le K\) form a compact
subset of \(\mathbb R\times X\). On either of the components of
\cref{thm_global_bifurcation} and \cref{thm_unilateral_global},
if the base flow is the only laminar member, then
\begin{equation}\label{eq_physical_compactness_alternative}
\sup_{\qty(Q,h)\text{ on the component}}\qty(M+V_I)=+\infty .
\end{equation}
\end{proposition}
\begin{proof}
The circulation bound controls \(Q\) since \(H_{\max}\le |p_0|M\). Equation \eqref{eq_physical_interface_pressure} now bounds the lower interface speed. The layer bounds \eqref{eq_physical_layer_speed_bounds} control \(h_p^{-1}\) and \(|h_q|/M\), so \cref{lem_gradient_to_closed_layer} gives uniform piecewise \(C^{3,\beta}\) bounds. Closedness of the constraint follows from \eqref{eq_physical_kinematics} and \(h_p\) has a uniform positive lower bound. Compact embedding into \(X\) gives compactness; the two continuation theorems then force a second laminar point on any component confined to such a set.
\end{proof}

\subsection{Interface speeds and layer separation}\label{subsec_interface_speed_separation}

A bound on \(Q\) alone keeps the interface a positive distance from
the bed. By considering the streamfunction on both sides of the interface, we obtain a local bound on the
interfacial speed. With both \(M\) and \(Q\) bounded, comparison and
energy estimates separate the interface from the free surface as well.

\begin{definition}\label{def_extremal_depths}
Write \(G_S=h(\cdot,0)\) and \(G_I=h(\cdot,p_I)\). With \(H_{\max}\)
as in \eqref{eq_physical_control_quantities}, set
\begin{equation}\label{eq_extremal_depths}
\begin{aligned}
b_1&=\min_{\mathbb T}G_I,& b_2&=\max_{\mathbb T}G_I,& H_{\min}&=\min_{\mathbb T}G_S,\\
a_1&=H_{\min}-b_1,& a_2&=H_{\max}-b_2,& a_*&=\max\{a_1,a_2\}.
\end{aligned}
\end{equation}
Both \(a_1\) and \(a_2\) are positive because \(G_S-G_I\) has
a positive minimum. For positive comparison depths \(a,b\), write
\begin{equation}\label{eq_comparison_layer_coefficients}
P(a)=\frac{\ell_+^2}{a^3}-\frac{\gamma_+^2a}{4},
\qquad
R(b)=\frac{\ell_-^2}{b^3}-\frac{\gamma_-^2b}{4}.
\end{equation}
\end{definition}

\begin{proposition}\label{prop_bernoulli_lower_depth}
Use the quantities in \cref{def_extremal_depths}. For \(t\in\mathbb R\)
write \(t^-=\max\{-t,0\}\) and put
\[
W_0=\max\left\{\sqrt Q+\gamma_+^-H_{\max},
                     \frac32\sqrt{\ell_+|\gamma_+|}\right\}.
\]
Every physical solution satisfies
\begin{equation}\label{eq_lower_depth_bernoulli_bound}
b_1\ge
\frac{\ell_-}{\gamma_-^-H_{\max}/2+
            \sqrt{rW_0^2+\delta(Q+2\gamma_+^-\ell_+)}}.
\end{equation}
Thus \(Q\le L\) gives a positive lower bound for \(b_1\) depending
only on \(L\) and the fixed data.
\end{proposition}
\begin{proof}
We bound both speed traces at a lowest point of the interface
in terms of \(Q\) and the extreme heights. Use \(z=y+D_0\), and
write \(s_{I,+}\) and \(s_{I,-}\) for the speed traces at that point.
Write \(s_{S,1}\) for the speed at a lowest point of the free surface.

For depth \(d\), flux \(\ell\) and vorticity \(\gamma\), use
\begin{equation}\label{eq_quadratic_stream_comparison}
w_{\ell,\gamma,d}(t)
=\ell(1-t/d)+\frac{\gamma}{2}t(d-t),\qquad 0\le t\le d.
\end{equation}
The profile satisfies \(w''=-\gamma\). Its downward derivatives
at \(0,d\) are \(\ell/d-\gamma d/2\) and \(\ell/d+\gamma d/2\),
respectively. In the upper layer choose \(d=a_1\) when
\(\gamma_+=0\), and otherwise set
\(d=\min\{a_1,\sqrt{\ell_+/|\gamma_+|}\}\).
The profile \(w_{\ell_+,\gamma_+,d}(z-b_1)\) decreases from
\(\ell_+\) to zero. On \(D_+\cap\{z<b_1+d\}\) it lies below
\(\psi_+\) by harmonic comparison. On the interface it is at
most \(\ell_+\), and it is zero on the horizontal cut, where
\(\psi_+\ge0\). The vertical quotient at the lowest interface
point gives \(s_{I,+}\le U_I(d)=\ell_+/d-\gamma_+d/2\).
When \(d=a_1\), a downward quotient at the lowest free-surface point
gives \(U_S(d)=\ell_+/d+\gamma_+d/2\le s_{S,1}\le\sqrt Q\).
The two choices of \(d\) give
\begin{equation*}
\begin{aligned}
U_I(d)&=U_S(d)-\gamma_+d\le\sqrt Q+\gamma_+^-H_{\max}
&&\text{for }d=a_1,\\
U_I(d)&=\sqrt{\ell_+|\gamma_+|}
        -\frac{\gamma_+}{2}\sqrt{\ell_+/|\gamma_+|}
        \le\frac32\sqrt{\ell_+|\gamma_+|}
&&\text{for }d<a_1 .
\end{aligned}
\end{equation*}
In \(0<z<b_1\) the profile \(w_{\ell_-,\gamma_-,b_1}(z)\) lies
below \(\psi_--\ell_+\), by equality at the bed and nonnegativity
of \(\psi_--\ell_+\) at the cut. The downward quotient at the
lowest interface point and \eqref{eq_physical_interface_pressure} give
\begin{equation*}
\begin{aligned}
\frac{\ell_-}{b_1}-\frac{\gamma_-^-H_{\max}}2
  &\le \frac{\ell_-}{b_1}+\frac{\gamma_-b_1}{2}\le s_{I,-},\\
s_{I,-}^2&=rs_{I,+}^2+\delta(Q-2\gamma_+\ell_+-2gb_1)\\
&\le rW_0^2+\delta(Q+2\gamma_+^-\ell_+).
\end{aligned}
\end{equation*}
These inequalities give \eqref{eq_lower_depth_bernoulli_bound}.
At a highest point of the free surface the Bernoulli equation gives
\(0<H_{\max}<Q/(2g)\), so the denominator is bounded when \(Q\le L\).
\end{proof}

\begin{proposition}\label{prop_interface_speed}
Let \(\qty({Q},{h})\in{\mathscr Z}\) be a physical solution as in \cref{prop_physical_speed_bounds}. Write
\({\ell}_{\max}=\max\{{\ell}_+,{\ell}_-\}\) and
\(A=\max\{|{\gamma}_+|,|{\gamma}_-|\}\), and assume
\(A>0\). At each \(x_0\in I\) define
\begin{equation}\label{eq_interface_radius}
{R}\qty(x_0)
=\min\{\operatorname{dist}\qty(x_0,S),\operatorname{dist}\qty(x_0,B),\pi\},
\end{equation}
with distances measured in the flat physical cylinder. There is a
constant \(C_*>0\) depending only on the dimension such that, for every
such solution and every \(x_0\in I\),
\begin{equation}\label{eq_interface_product_bound}
|\nabla{\psi}_+\qty(x_0)|\,|\nabla{\psi}_-\qty(x_0)|
\le C_*\max\left\{\frac{{\ell}_{\max}^2}{{R}\qty(x_0)^2},
A{\ell}_{\max}\right\}
\end{equation}
and
\begin{equation}\label{eq_interface_upper_speed_bound}
|\nabla{\psi}_+\qty(x_0)|^2
\le \frac{2{\delta}\max\{{\gamma}_+,0\}{\ell}_+}{{r}}
+\frac{C_*}{\sqrt{{r}}}
\max\left\{\frac{{\ell}_{\max}^2}{{R}\qty(x_0)^2},
A{\ell}_{\max}\right\}.
\end{equation}
The laminar family of \cref{thm_laminar_family} with the same fixed
data attains the inverse-square scale in \eqref{eq_interface_product_bound}.
For every choice of interface point \(x_0\) along that family,
\begin{equation}\label{eq_interface_product_sharpness}
\lim_{{\lambda}\to\infty}
{R}\qty(x_0)^2
|\nabla{\psi}_+\qty(x_0)|\,|\nabla{\psi}_-\qty(x_0)|
=\min\{{\ell}_+,{\ell}_-\}^2 .
\end{equation}
\end{proposition}
\begin{proof}
Set \(u_-=(\psi-\ell_+)_+\), \(u_+=(\ell_+-\psi)_+\). Each lies between zero and \(\ell_i\), is supported on its own phase, and extends continuously by zero in \(H^1_{\rm loc}\). Its interior Laplacian is \(-\gamma_-\) or \(\gamma_+\), respectively. Integration by parts adds the nonnegative interface measure \(|\nabla\psi_i|\,ds\), hence \(\Delta u_i\ge-A\). Set \(R_0=\min(R(x_0),\sqrt{\ell_{\max}/A})\), \(v_i(z)=u_i(x_0+R_0z)/\ell_{\max}\). We can lift the disk to the periodic fluid domain since it meets neither outer boundary. On \(B_1\), \(0\le v_i\le1\), \(v_-v_+=0\), \(\Delta v_i\ge-1\). Testing with \(\eta^2v_i\) for a fixed cutoff equal to one on \(B_{1/2}\) and absorbing the mixed term bounds \(\int_{B_{1/2}}|\nabla v_i|^2\le C\).

For \(w_i(z)=v_i(z/2)\) the unit-disk energies are bounded by planar scaling. The two-phase estimate \cite[Theorem~3.1]{Velichkov2014} (see also \cite{AltCaffarelliFriedman1984,CaffarelliJerisonKenig2002}) gives \(\prod_i(s^{-2}\int_{B_s}|\nabla w_i|^2)\le C\) for \(0<s<1\). The \(C^1\) interface bisects the disk asymptotically and each gradient has a one-sided limit. Thus each factor tends to \((\pi/2)(R_0/(2\ell_{\max}))^2|\nabla\psi_i(x_0)|^2\). Taking square roots, we obtain \eqref{eq_interface_product_bound}. Write \(X=|\nabla\psi_+|^2\), \(\chi=2\delta\max(\gamma_+,0)\ell_+\). The surface law gives \(Q-2gh_I>0\); hence \eqref{eq_physical_interface_pressure} gives \(|\nabla\psi_-|^2\ge rX-\chi\). For \(X>\chi/r\), \(r(X-\chi/r)^2\le X(rX-\chi)\); the product bound gives the asserted upper speed estimate, which is immediate for \(X\le\chi/r\). Finally the laminar speeds are uniformly asymptotic to \(\sqrt\lambda\), and \(d_i\sqrt\lambda\to\ell_i\). Thus \(R=\min(d_+,d_-,\pi)\) gives the stated sharpness limit.
\end{proof}

\begin{proposition}\label{prop_nonstagnant_gap_control}
Use the physical notation of \cref{def_physical_controls}.
For \(K,L>0\), the physical solutions with \(M\le K\) and \(Q\le L\)
form a compact subset of \(\mathbb R\times X\). There are
\(d_0>0\) and \(C<\infty\), depending only on \(K,L\) and the fixed
data, such that these solutions satisfy
\begin{equation}\label{eq_nonstagnant_gap_control}
{d_*}\ge d_0,\qquad
V_I+{B_I}\le C .
\end{equation}
\end{proposition}
\begin{proof}
We rescale the streamfunction in the lower layer to obtain two
distributional inequalities, use vertical comparisons to bound the
distance of both graphs from the bed, and compare disk averages to
separate the graphs from each other. The separation controls the
interface speed and slope and gives compactness.

\smallskip\noindent\textit{Step 1. Distributional inequalities for the rescaled streamfunctions.}\enspace
For \(M\le K\) and \(Q\le L\), we obtain
\eqref{eq_comparison_heads_laplacians} by choosing one scale with
nonpositive interface flux density and one with nonnegative density.
Set \(z=y+D_0\), \(G_S=h\qty(\cdot,0)\), \(G_I=h\qty(\cdot,p_I)\) and
\(T=|p_0|K\). By \cref{prop_physical_speed_bounds},
 \(0<G_I<G_S\le T,\qquad |G_S'|\le K\sqrt L .\) 
Write \(s_i=|\nabla\psi_i|\) at the interface, so \(s_i\ge K^{-1}\).
For the interface pressure relation \eqref{eq_physical_interface_pressure},
use \(Q-2gG_I\ge Q-2gG_S\ge0\) and \(Q\le L\) to obtain
\begin{equation*}
\begin{aligned}
r s_+^2-b_0&\le s_-^2\le r s_+^2+b_1,&
b_0&=2\delta\max\{\gamma_+,0\}\ell_+,&
b_1&=\delta\qty(L+2|\gamma_+|\ell_+).
\end{aligned}
\end{equation*}
Absorbing \(b_0\) with \(s_-^2\ge K^{-2}\) and \(b_1\) with
\(s_+^2\ge K^{-2}\) gives
\begin{equation*}
R_0\le s_-/s_+\le R_1,\qquad
R_0=\sqrt{\frac r{1+b_0K^2}},\qquad R_1=\sqrt{r+b_1K^2}.
\end{equation*}
We set \(a_0=\min\{1,1/R_1\}\), \(a_1=\max\{1,1/R_0\}\) and
\(A=\max\{|\gamma_+|,a_1|\gamma_-|\}\). The fixed vorticities are
distinct, so \(A>0\). For \(a\in\{a_0,a_1\}\), define on the
half-cylinder \(z>0\)
\begin{equation*}
\chi_a\qty(x,z)=
\begin{cases}
0,&z\ge G_S\qty(x),\\
\psi\qty(x,z),&G_I\qty(x)<z<G_S\qty(x),\\
\ell_++a\qty(\psi\qty(x,z)-\ell_+),&0<z\le G_I\qty(x).
\end{cases}
\end{equation*}
The streamfunction in these coordinates is evaluated at physical
height \(y=z-D_0\). Its boundary values give continuity and local
\(H^1\) regularity, and \(0\le\chi_{a_1}\le\qty(a_1/a_0)\chi_{a_0}\).
The ordinary Laplacian is \(-\gamma_+\) in the upper layer,
\(-a\gamma_-\) in the lower layer and zero above the surface.
With the normal pointing upward at the interface,
the two normal derivatives of \(\psi\) are \(-s_i\).
Thus the distributional Laplacian of \(\chi_a\) has interface
density \(as_--s_+\) and surface density \(s_+|_S\in\qty[0,\sqrt L]\).
The bounds on the speed ratio give
\begin{equation}\label{eq_comparison_heads_laplacians}
\Delta\chi_{a_1}\ge -A,\qquad
\Delta\chi_{a_0}\le A+\sqrt L\,\mathcal H^1\lfloor S.
\end{equation}

\smallskip\noindent\textit{Step 2. Separation from the bed.}\enspace By \cref{prop_bernoulli_lower_depth}, the minimum interface height above the bed has a positive lower bound depending only on \(L\) and the data. The free surface lies above the interface pointwise and inherits this lower bound. Denote the common positive lower bound by \(c_I=c_S\). All disks used below can therefore be chosen away from the bed.

\smallskip\noindent\textit{Step 3. Separation of the interface and free surface.}\enspace
We argue by contradiction. A sufficiently small fixed distance \(R\) between
points of \(I\) and \(S\) would force incompatible lower and upper
bounds on a disk mean. We set \(J=\sqrt{1+K^2L}\) and choose \(R>0\) so
small that
\begin{equation*}
\begin{aligned}
R&<\min\{\pi/4,c_S/4\},\qquad AR^2/4\le\ell_+/4,\\
4\qty(a_1/a_0)\qty(AR^2+2\sqrt L\,JR/\pi)&\le\ell_+/4.
\end{aligned}
\end{equation*}
If \(X_I\in I\) and \(X_S\in S\) were at distance at most \(R\),
their nearest planar lifts would satisfy
\(B_R\qty(X_I)\subset B_{2R}\qty(X_S)\subset\{z>0\}\), with the larger disk
injecting into the cylinder. By \eqref{eq_comparison_heads_laplacians},
\(\chi_{a_1}+A|X-X_I|^2/4\) is subharmonic. The disk mean-value
inequality and \(\chi_{a_1}\qty(X_I)=\ell_+\) give
 \(\fint_{B_R\qty(X_I)}\chi_{a_1}\ge\ell_+-AR^2/4\ge3\ell_+/4.\) 
For the upper bound let \(m\qty(t)\) be the circular mean of \(\chi_{a_0}\)
on \(\partial B_t\qty(X_S)\) and let \(\mu=\Delta\chi_{a_0}\) as a signed
measure in this disk. The slope bound gives
\(\mathcal H^1\qty(S\cap B_t\qty(X_S))\le2Jt\). For \(0<s<t<2R\),
\begin{equation*}
\begin{aligned}
m\qty(t)-m\qty(s)
&=\frac1{2\pi}\int_s^t\frac{\mu\qty(B_\tau\qty(X_S))}{\tau}\,d\tau\\
&\le \frac A4\qty(t^2-s^2)+\frac{\sqrt L\,J}{\pi}\qty(t-s).
\end{aligned}
\end{equation*}
The identity follows from the divergence theorem for mollifications.
The means converge by continuity, and the disk masses converge for
almost every radius since the corresponding circles have zero total
variation mass. On \(\qty[s,t]\) the masses are bounded by the total
variation on a slightly larger disk. Taking \(s\downarrow0\) uses
\(\chi_{a_0}\qty(X_S)=0\). Averaging the resulting bound in \(t\) and then
using the disk inclusion gives
\begin{equation*}
\begin{aligned}
\fint_{B_{2R}\qty(X_S)}\chi_{a_0}&\le AR^2+2\sqrt L\,JR/\pi,\\
\fint_{B_R\qty(X_I)}\chi_{a_1}
&\le4\qty(a_1/a_0)\fint_{B_{2R}\qty(X_S)}\chi_{a_0}\le\ell_+/4 .
\end{aligned}
\end{equation*}
This is incompatible with the lower bound already obtained on the disk centered at \(X_I\).
The distance from \(I\) to \(S\) is therefore at least \(R\), and
\({d_*}\) is at least \(\min\{c_I,R\}\).

\smallskip\noindent\textit{Step 4. Interface bounds and compactness.}\enspace
The lower bound for \({d_*}\) gives the remaining
bound in \eqref{eq_nonstagnant_gap_control} and places the solutions
with \(M\le K\) and \(Q\le L\) in a compact sublevel.
The interface speed estimate in \cref{prop_interface_speed} now
bounds \(V_I\), and \({B_I}\le M V_I\) by the
kinematic identity \eqref{eq_physical_kinematics}. The solution set
lies in a compact \(M+V_I\) sublevel by
\cref{prop_physical_compactness}. Within that set the inequalities
\(M\le K\) and \(Q\le L\) are closed, proving compactness.
\end{proof}

\subsection{Depth and the Bernoulli bound}\label{subsec_shallow_layer_rigidity}

At bounded \(M\), large \(Q\) forces small total depth. We use a weighted trace estimate to exclude that regime for every nonlaminar solution.

\begin{proposition}\label{prop_fast_flow_thinning}
Use the physical quantities of \cref{prop_physical_speed_bounds}.
For \(K>0\), put
\begin{equation}\label{eq_fast_flow_threshold}
\begin{aligned}
T&=|p_0|K,& A_\pm&=\max\{\gamma_\pm,0\},&
b_0&=2gT+2|\gamma_+|\ell_+,\\
{Q_0}\qty(K)&=1+2b_0+2T^2\qty(A_-^2/\delta+A_+^2).
\end{aligned}
\end{equation}
Write \(H_I=\max_{q\in\mathbb T}h\qty(q,p_I)\). For every solution with
\(M\le K\) and \(Q\ge {Q_0}\qty(K)\),
\begin{equation}\label{eq_fast_flow_depth_bounds}
\begin{aligned}
H_I&\le\frac{2\sqrt{2}\ell_-}{\sqrt{\delta Q}},&
H_{\max}-H_I&\le\frac{2\sqrt{2}\ell_+}{\sqrt Q},\\
H_{\max}&\le\frac{2\sqrt{2}\qty(\ell_++\ell_-/\sqrt{\delta})}{\sqrt Q}.
\end{aligned}
\end{equation}
In particular, the two pointwise vertical thicknesses and the total
depth converge uniformly to zero along any sequence with bounded
\(M\) and \(Q\to+\infty\).
\end{proposition}
\begin{proof}
We compare the streamfunction with one-dimensional profiles at the
highest point of the interface and at the highest point of the surface.
Work with \(z=y+D_0\), so \(0<H_I<H_{\max}\le T\). Write \(s_i=|\nabla\psi_i|\).
The Bernoulli and pressure identities in \cref{prop_physical_speed_bounds}
give
\begin{equation*}
\begin{aligned}
s_+^2|_S&=Q-2gz\ge Q-b_0\ge Q/2,\\
s_-^2|_I&=r s_+^2|_I+\delta\qty(Q+2\gamma_+p_I-2gz)
\ge\delta\qty(Q-b_0)\ge\delta Q/2 .
\end{aligned}
\end{equation*}
In the lower layer \(w=\psi-\ell_+\) is bounded above by
\begin{equation*}
b_-\qty(z)=\ell_-\qty(1-\frac z{H_I})+\frac{A_-}2z\qty(H_I-z),
\qquad 0\le z\le H_I .
\end{equation*}
Indeed \(w=b_-\) at the bed and \(w=0\le b_-\) on \(I\), while
\(\Delta\qty(w-b_-)=A_--\gamma_-\ge0\). At any highest point of \(I\),
the interface is horizontal, and both profiles vanish. A downward
difference quotient there bounds the normal speed:
\begin{equation*}
\sqrt{\delta Q/2}\le s_-\le\frac{\ell_-}{H_I}+\frac{A_-H_I}{2},
\qquad
\frac{\ell_-}{H_I}\ge\frac12\sqrt{\delta Q/2}.
\end{equation*}
The last inequality uses \(H_I\le T\) and \(Q\ge2A_-^2T^2/\delta\).

In the portion of the upper layer with \(z>H_I\), put \(L=H_{\max}-H_I\)
and compare \(\psi\) with
\begin{equation*}
b_+\qty(z)=\ell_+\qty(1-\frac{z-H_I}{L})
+\frac{A_+}2\qty(z-H_I)\qty(H_{\max}-z),
\qquad H_I\le z\le H_{\max}.
\end{equation*}
The difference \(\psi-b_+\) is subharmonic. It is nonpositive on the
horizontal cut \(z=H_I\) because \(\psi\le\ell_+\), and on the surface
because \(b_+\ge0\). The maximum principle applies on every component
of the cut layer in the periodic cylinder. At a highest point of the surface the downward
quotient therefore gives
\begin{equation*}
\sqrt{Q/2}\le s_+\le\frac{\ell_+}{L}+\frac{A_+L}{2},
\qquad
\frac{\ell_+}{L}\ge\frac12\sqrt{Q/2}.
\end{equation*}
Here \(L\le T\) and \(Q\ge2A_+^2T^2\). Combining the two bounds gives
\eqref{eq_fast_flow_depth_bounds} because \(H_{\max}=H_I+L\).
Each pointwise layer thickness is bounded above by \(H_{\max}\).
\end{proof}

\begin{proposition}\label{prop_nodal_fast_flow_bound}
Every nonlaminar physical solution satisfies
\begin{equation}\label{eq_nodal_depth_and_q}
 gM^2H_{\max}>1,\qquad M>(g|p_0|)^{-1/3}.
\end{equation}
\end{proposition}
\begin{proof}
Tangentially differentiating the equation gives the homogeneous conormal problem for \(W=-h_q\) with coefficients \eqref{eq_nodal_tangential_system}. By periodicity, the regularity already established and the zero bed trace, we may take \(W\) as the test function. Completing the square in the coefficient matrix gives
\[
 E=\sum_i\rho_i\int_{\Omega_i}\left[
 \frac{1}{h_p}\left(W_q-\frac{h_q}{h_p}W_p\right)^2+
 \frac{W_p^2}{h_p^3}\right]
 =g\rho_+\int_{\mathbb T}W_S^2+g(\rho_--\rho_+)\int_{\mathbb T}W_I^2 .
\]
For \(t=0,p_I\), the common trace across the interface and the zero bed trace permit a single vertical integration. By Cauchy--Schwarz,
\[
 W(q,t)^2\le \left(\int_{p_0}^t h_p^3\right)
 \left(\int_{p_0}^t W_p^2h_p^{-3}\right)
 \le M^2h(q,t)\int_{p_0}^t W_p^2h_p^{-3}.
\]
Both heights are bounded above by \(H_{\max}\). The contribution of the surface trace carries weight \(\rho_+\) on both layers; the interface trace adds precisely \(\rho_--\rho_+\) on the lower layer. Thus \(E\le gM^2H_{\max}\sum_i\rho_i\int_{\Omega_i}W_p^2 h_p^{-3}\le gM^2H_{\max}E\). For a nonlaminar wave \(W\ne0\), and \(E>0\) by ellipticity and the zero bed trace. To prove that the inequality is strict, suppose \(gM^2H_{\max}=1\). Since \(\max_q h(q,p_I)<H_{\max}\) and the interface weight is positive, equality in the interface estimate forces \(\int_{\Omega_-} W_p^2h_p^{-3}=0\). It then follows from the bed condition that \(W=0\) in the lower layer. Integrating from \(p_I\) to the surface now bounds the top trace with the smaller factor \(M^2(H_{\max}-\min_q h(q,p_I))\). Since the interface stays strictly above the bed for this wave, the energy identity gives \(E<E\), a contradiction. This strict inequality is pointwise in the solution; the argument supplies no uniform margin. Now \(H_{\max}\le |p_0|M\) gives the second.
\end{proof}
\begin{corollary}\label{cor_nonlaminar_bernoulli_bound}
If \(M\le K\), every nonlaminar solution satisfies
 \(Q\le Q_{\rm c}(K):=\max\{Q_0(K),C^2g^2K^4\},\qquad C=2\sqrt2(\ell_++\ell_-/\sqrt\delta).\) 
\end{corollary}
\begin{proof}
When \(Q\ge Q_0(K)\), \eqref{eq_fast_flow_depth_bounds} gives \(H_{\max}\le C/\sqrt Q\). We apply \eqref{eq_nodal_depth_and_q}. The complementary range is covered by the first term in the maximum.
\end{proof}

\subsection{Global alternatives and horizontal stagnation}\label{subsec_physical_dichotomy}

\begin{proposition}\label{thm_physical_global_alternatives}
For every \(K>0\), the set \(\{z\in\mathscr S:M(z)\le K\}\) is compact in \(\mathbb R\times X\). On this set the interfacial speed, the interface slope and the maximum wave-frame speed remain bounded, and the boundary separation is bounded below by a positive number.
\end{proposition}
\begin{proof}
Approximate an arbitrary point of this sublevel by nonlaminar solutions. Continuity of \(M\) eventually gives \(M\le K+1\) on the approximants; \cref{cor_nonlaminar_bernoulli_bound} gives \(Q\le Q_{\rm c}(K+1)\) at the limit. The sublevel in question is relatively closed in the compact physical zero set of \cref{prop_nonstagnant_gap_control} for these bounds. Combining that proposition with \cref{prop_physical_speed_bounds} gives the stated uniform bounds.
\end{proof}
\begin{proof}[Proof of \cref{thm_bilateral_stagnation} and \cref{thm_oriented_stagnation}]
Both components are relatively closed in the physical open set, so the preceding proposition applies to their sublevels. Nodal strictness follows from \cref{thm_nodal_positivity}; \cref{prop_laminar_return} gives the spectral condition on a second laminar point in each case. If the base flow is the only laminar point on the component, the continuation alternatives in \cref{thm_global_bifurcation} and \cref{thm_unilateral_global} exclude compactness. Consequently \(M\) is unbounded. Along a sequence with \(M\to\infty\), the identity \(c-u_i=1/h_p\) gives \(\min(c-u)\to0\), as required. Uniqueness excludes the oriented return, and \cref{cor_fold} gives the structure of the bilateral component.
\end{proof}

The minimizing point in \(\min(c-u)\) may move with the wave. This conclusion leaves open both the existence of a limiting wave with vanishing full relative velocity and the location at which degeneration occurs. When there is more than one root, the component may contain both a second laminar point and an unbounded sequence.

\subsection{Pressure at the interface}\label{subsec_interface_pressure}

We conclude with a separate criterion for the sign of the pressure at the interface. For a single homogeneous fluid, V\u{a}rv\u{a}ruc\u{a} \cite[Theorem~2.1]{Varvaruca2008} proved a pressure sign estimate for nonpositive vorticity. The comparison in \cref{prop_interface_pressure_positivity} uses the common pressure trace and a normalization that matches the flux of the pressure quotient across the density jump.

\begin{proposition}\label{prop_interface_pressure_positivity}
For a physical solution, write the shifted height and positive
horizontal relative speed as
 \(z=y+D_0,\qquad w_i=c-u_i=-\psi_{i,z}>0 .\) 
Define the reference head and its minimum by
\begin{equation}\label{eq_pressure_reference_head}
H_{\rm P}(z)=Q+2\gamma_+p_I-2gz,\qquad
H_{\rm P,*}=Q+2\gamma_+p_I-2gH_{\max}.
\end{equation}
If
\begin{equation}\label{eq_pressure_sign_conditions}
H_{\rm P,*}>0,\qquad
\gamma_i w_i\le g\quad\text{on }\overline D_i,\quad i=\pm,
\end{equation}
then \(P_i\ge0\) on both closed layers. If \(\gamma_+,\gamma_-\le0\),
these conditions hold and the pressure is strictly positive on the
entire closed fluid domain away from the free surface. In particular,
\begin{equation}\label{eq_pressure_strict_interface_speed}
|\nabla\psi_+|^2<Q+2\gamma_+p_I-2gh(q,p_I)
\qquad\text{on }I.
\end{equation}
\end{proposition}
\begin{proof}
We introduce a quotient by dividing the pressure by a positive function of depth. The normal
derivatives of this quotient satisfy a homogeneous transmission law, which excludes
a negative minimum at the interface as well as in the layers.

Set \(\Pi_i=P_i/\rho_i\) and \(s_i=|\nabla\psi_i|\).
The Bernoulli functions in the change of variables \(p=-\psi\) give
\(s_i^2=Q+2b_i(p)-2gz-2\Pi_i\) with \(p=-\psi_i\) and \(b_i'=\gamma_i\).
Write \(F_i=(-\psi_{i,z},\psi_{i,x})=s_it_i\), with
\(t_i=(\cos\vartheta_i,\sin\vartheta_i)\),
\(\vartheta_i=\arctan(\psi_{i,x}/w_i)\), and
\(n_i=(-\sin\vartheta_i,\cos\vartheta_i)\).
The rotation \(J(a,b)=(b,-a)\) sends \(n_i\) to \(t_i\). Using
\(\operatorname{div}F_i=0\) in the Euler equation gives
\begin{equation*}
\begin{aligned}
\nabla(\Pi_i+gz)
&=-(F_i\cdot\nabla)F_i
=-s_i(\partial_{t_i}s_i)t_i-s_i^2(\partial_{t_i}\vartheta_i)n_i\\
&=s_i^2\bigl((\partial_{n_i}\vartheta_i)t_i
-(\partial_{t_i}\vartheta_i)n_i\bigr)=s_i^2J\nabla\vartheta_i,\\
\operatorname{div}\bigl(s_i^{-2}\nabla(\Pi_i+gz)\bigr)&=0 .
\end{aligned}
\end{equation*}
We set \(R_{\rm P}=P/H_{\rm P}\) on the closed fluid domain.
The denominator is positive by \eqref{eq_pressure_sign_conditions}.
The quotient is continuous across \(I\), is \(C^{1,\alpha}\) on each
closed layer, and is smooth in each open layer. Use
\(\Pi_i=H_{\rm P}R_{\rm P}/\rho_i\) and \(p_z=w_i\). Differentiating the
Bernoulli expression gives
\begin{equation*}
\partial_z s_i^2=2\gamma_iw_i-2g
+\frac{4gR_{\rm P}}{\rho_i}-\frac{2H_{\rm P}}{\rho_i}\partial_z R_{\rm P}.
\end{equation*}
Substitution in the divergence identity yields
\begin{equation}\label{eq_pressure_ratio_layer}
\begin{aligned}
\Delta R_{\rm P}+\mathbf b_i\cdot\nabla R_{\rm P}
&=\frac{2g(\rho_i-2R_{\rm P})}{H_{\rm P}s_i^2}
\left(\gamma_iw_i-g\left(1-\frac{2R_{\rm P}}{\rho_i}\right)\right),\\
\mathbf b_i
&=-\nabla\log s_i^2-\frac{4g}{H_{\rm P}}e_z
+\frac{2g(\rho_i-2R_{\rm P})}{\rho_i s_i^2}e_z.
\end{aligned}
\end{equation}
On a fixed wave, \(s_i\ge w_i>0\) and \(H_{\rm P}\ge H_{\rm P,*}>0\)
on the closed layers. Thus the drift is bounded.
Where \(R_{\rm P}<0\), \eqref{eq_pressure_sign_conditions}
makes the right-hand side strictly negative.

The common tangent \(t\) and upward normal \(n\) at the interface
give the same curvature \(\partial_t\vartheta_i\) on both sides.
The pressure gradient identity and the common Bernoulli offset at
\(I\) yield
\begin{equation*}
\begin{aligned}
\rho_i^{-1}\partial_n P_i&=-g n_z-s_i^2\partial_t\vartheta_i,&
s_i^2&=H_{\rm P}-2P/\rho_i,\\
s_-^2-s_+^2&=2(1/\rho_+-1/\rho_-)P,&
\rho_-s_-^2-\rho_+s_+^2&=(\rho_--\rho_+)H_{\rm P}.
\end{aligned}
\end{equation*}
Use \(a_i=(\rho_i s_i^2)^{-1}\) and
\(d=2g n_z(1/\rho_+-1/\rho_-)/(s_+^2s_-^2)\) on \(I\).
The curvature terms cancel, and \(\partial_n H_{\rm P}=-2g n_z\).
Consequently,
\begin{equation}\label{eq_pressure_ratio_interface}
\begin{aligned}
a_-\partial_n P_- - a_+\partial_n P_+&=dP,\\
a_-\partial_n H_{\rm P}-a_+\partial_n H_{\rm P}&=dH_{\rm P},\\
a_-\partial_n R_{{\rm P},-}-a_+\partial_n R_{{\rm P},+}&=0 .
\end{aligned}
\end{equation}
Here we take all normal derivatives with respect to the same upward normal.
The last equality follows from the quotient rule, since the two traces have common \(P\) and \(H_{\rm P}\).

A negative global minimum of \(R_{\rm P}\) cannot occur in an open layer
by \eqref{eq_pressure_ratio_layer}, or on the free surface where
\(R_{\rm P}=0\). At the bed, \(\vartheta_-=0\) and its tangential derivative
vanishes. The gradient formula gives \(\partial_z P_-=-\rho_-g\), hence
\begin{equation*}
\partial_z R_{{\rm P},-}=\frac{g(2R_{\rm P}-\rho_-)}{H_{\rm P}}<0
\end{equation*}
at a negative bed minimum, contradicting its inward derivative.
At a negative interface minimum there is a neighborhood where
\(R_{\rm P}<0\) in both layers. The strict supersolution inequality
and the interior ball condition on both sides give opposite Hopf
derivatives \(\partial_n R_{{\rm P},-}<0<\partial_n R_{{\rm P},+}\).
This contradicts \eqref{eq_pressure_ratio_interface} because \(a_i>0\).
It follows that \(P_i\ge0\).

When \(\gamma_i\le0\), \(2\gamma_+p_I\ge0\) and the surface identity
gives \(Q-2gH_{\max}>0\). Also \(\gamma_i w_i\le0\).
If \(R_{\rm P}=0\) at a minimum away from the surface, there is a
neighborhood with \(0\le R_{\rm P}<\min(\rho_+,\rho_-)/4\).
The right-hand side of \eqref{eq_pressure_ratio_layer} is strictly
negative there, and the same interior, bed and interface arguments
exclude the zero minimum. Equation \eqref{eq_pressure_strict_interface_speed}
now follows from \(s_+^2=H_{\rm P}-2P/\rho_+\) on \(I\).
\end{proof}
\begin{corollary}\label{cor_pressure_gap}
If \(\gamma_+,\gamma_-\le0\) and \(Q\le L\), then \(s_i=|\nabla\psi_i|\le C(L)\) throughout the closed layers. In particular
 \[\operatorname{dist}(I,S)\ge\ell_+/C(L),\qquad \operatorname{dist}(I,B)\ge\ell_-/C(L).\] 
\end{corollary}
\begin{proof}
At the interface \cref{prop_interface_pressure_positivity} gives \(P_I>0\).
The Bernoulli laws imply \(s_i^2|_I=Q+2\gamma_+p_I-2gh_I-2P_I/\rho_i
 \le L+2|\gamma_+|\ell_+\). At the surface \(s_+^2\le L\), and \(H_{\max}\le L/(2g)\). The upper maximum principle and lower reflection estimate of \cref{prop_physical_speed_bounds} therefore bound \(s_i\) throughout. For each pair of adjacent boundary graphs, choose points realizing their distance in the periodic cylinder. The interior of a shortest joining segment lies in the intervening layer: an additional intersection would yield a shorter boundary-to-boundary segment. Integrating \(\nabla\psi_i\) along it gives the stated bounds. If \(M\) is bounded as well, \(h_p^{-1}\le C(L)\) and \(|h_q|\le MC(L)\); \cref{lem_gradient_to_closed_layer} then gives the closed-layer estimates directly.
\end{proof}

\end{document}